\documentclass[a4paper,final]{elsarticle}%
\usepackage{amssymb}
\usepackage{amsfonts}
\usepackage{graphicx}
\usepackage{amsmath}
\usepackage{amssymb}
\usepackage{epsfig}%
\providecommand{\U}[1]{\protect\rule{.1in}{.1in}}
\providecommand{\U}[1]{\protect\rule{.1in}{.1in}}
\providecommand{\U}[1]{\protect\rule{.1in}{.1in}}
\newtheorem{defn}{Definition}[section]
\newtheorem{thm}[defn]{Theorem}
\newtheorem{exmp}[defn]{Example}

\newtheorem{rmk}[defn]{Remark}

\newtheorem{prop}[defn]{Proposition}
\newtheorem{cor}[defn]{Corollary}

\newtheorem{notation}[defn]{Notation}
\newenvironment{pf}[1][Proof]{\textbf{#1.} }{\ \rule{0.5em}{0.5em}}
\newenvironment{akm}[1][Acknowledgement]{\textbf{#1.} }{}

\journal{ArXiv}
\begin{document}
\begin{frontmatter}
\title{Products of $2\times2$ matrices related to non autonomous Fibonacci difference equations}
\author[A]{Rafael Lu\'{\i}s\footnote{Email: rafel.luis.madeira@gmail.com}}
\author[A]{Henrique M. Oliveira\footnote{Email: holiv@math.ist.utl.pt}}
\address[A]{
Center for Mathematical Analysis, Geometry and Dynamical Systems\\
Department of Mathematics - Technical Institute of Lisbon - University of Lisbon\\
Av. Rovisco Pais, 1049-001 Lisboa  - Portugal}
%\email{rafel.luis.madeira@gmail.com holiv@math.ist.utl.pt}
\begin{abstract}
A technique to compute arbitrary products of a class of Fibonacci  $2\times2$ square matrices is proved
in this work.
General explicit solutions for non autonomous Fibonacci difference equations are obtained from these products.
In the periodic non autonomous Fibonacci difference equations the monodromy matrix, the Floquet multipliers and the
Binet's formulas are obtained.
In the periodic case explicit solutions are obtained and the solutions are analyzed.
\end{abstract}
\end{frontmatter}

%\linenumbers


\section{Introduction}

\subsection{Historical background}

Difference equations can model effectively almost any physical and artificial
phenomena \cite{ELA2003}. One of the earliest recurrences, or in other words,
a difference equation, giving us the Fibonacci sequence, was introduced in
1202 in the old \textquotedblleft Liber Abaci" \cite{Pis2002}, a book about
the abacus, by the famous Italian mathematician Leonardo Pisano better known
as Fibonacci. The first numbers of this sequence are%
\[
0,1,1,2,3,5,8,13,21,\ldots.
\]
This example is related to ecology but one may find many applications of this
sequence of numbers in various branches of science like pure and applied
mathematics, in biology or in phyllotaxies, among many others.

If $x_{n}$ represents a number in this sequence, with $n=0,1,2,\ldots$, then
the Fibonacci numbers satisfy the recurrence%
\begin{equation}
q_{n+2}=q_{n+1}+q_{n},\text{ }q_{0}=0,\text{ }q_{1}=1. \label{Fib_DE}%
\end{equation}
Obviously, we can look the above recurrence as a difference equation. The
general solution of equation (\ref{Fib_DE}) is given by the Binet formula
\[
q_{n}=\frac{1}{\sqrt{5}}\left(  \phi^{n}-\left(  1-\phi\right)  ^{n}\right)
,
\]
where $\phi=\frac{1+\sqrt{5}}{2}$ is the well known golden ratio.

The equation (\ref{Fib_DE}) can be written in matrix notation%
\begin{equation}
X_{n+1}=f(X_{n}), \label{FIB_eq2}%
\end{equation}
where $f(X_{n})=AX_{n}$, $A=\left[
\begin{array}
[c]{cc}%
1 & 1\\
1 & 0
\end{array}
\right]  $, $X_{n}=\left[
\begin{array}
[c]{c}%
q_{n+1}\\
q_{n}%
\end{array}
\right]  $ and $X_{n+1}=\left[
\begin{array}
[c]{c}%
q_{n+2}\\
q_{n+1}%
\end{array}
\right]  $. It is easy to find the solution of equation (\ref{FIB_eq2}) via
the powers of the matrix $A$, i.e., $X_{n}=A^{n}X_{0}$ and using the
eigenvalues of $A$, $\phi$ and $1-\phi$.

Recently, Edson and Yayenie \cite{EY2010} studied a non-autonomous
generalization of the Fibonacci recurrence. They considered the equation
\begin{equation}
x_{n+2}=a_{n}x_{n+1}+x_{n}, \label{naFre}%
\end{equation}
with initial conditions $x_{0}$ and $x_{1}$ in which $a_{n}$ is a $2$-periodic
non-zero sequence of non-negative integer numbers, i.e., $a_{n}=a_{0}$ if $n$
is even and $a_{n}=a_{1}$ if $n$ is odd. They found a Binet Formula for this
equation using generating functions.

Now, let $\left\{  a_{n}\right\}  $, $n=0,1,2,\ldots$ be a $k$-periodic
non-zero sequence of non-negative integer numbers with fixed $k>2$, i.e.,
$a_{n}=a_{n+k}$, for all $n=0,1,2,\ldots$. It was stated in \cite{EY2010} that
to find a Binet formula for equation (\ref{naFre}), for all $n=0,1,2,\ldots$,
was an open problem.

Later on, Lewis joined the previous authors in another nice paper
\cite{EY2011} on the same subject, where the $k$-periodic equation was treated
using again generating functions. In that paper is developed an elegant new
technique to obtain the solutions, which are not presented explicitly for the
general case $\left\{  a_{n}\right\}  $, $n=0,1,2,\ldots$.

In a series of papers \cite{Mallik1997, Mallik1998, Mallik2000} Mallik studied
the solutions of linear difference equations with variable coefficients. Using
classical techniques of iteration, he was able to find the solutions of
certain special linear difference equations. In \cite{Mallik1997} he studied a
second order equation and he was able to write expressions for the solutions.
In \cite{Mallik1998} the author presented the solution of a linear difference
equation of unbounded order. As special cases, the solutions of nonhomogeneous
and homogeneous linear difference equations of order $n$ with variable
coefficients were obtained. From these solutions, Mallik was able to get the
expressions for the product of companion matrices, and the power of a
companion matrix. The closed form solution of the $n$-th order difference
equation $(n\geq3)$ is presented in \cite{Mallik2000} using some combinatorial
properties in the indices of the coefficients in an indirect manner. The
results of Mallik, being very interesting and original, were new approaches to
the classic problem of solving linear difference equations with variable coefficients.

In \cite{Jesus} de Jesus and Petronilho established algebraic conditions for
the existence of a polynomial mapping using a monic orthogonal polynomial
sequence (OPS). In particular in section 5.1 of this paper is introduced a
method to study such sequences with periodic coefficients that are related
with difference equations with periodic coefficients. These results were
developed in a recent work \cite{Petronilho} by Petronilho to compute the
solutions of the periodic case via orthogonal polynomials and a determinant of
a tridiagonal matrix associated with the dynamics generated by equation
\ref{naFre}. Again, an explicit expression of the solution for the general
case $\left\{  a_{n}\right\}  $, $n=0,1,2,\ldots$,\ is not presented. We point
out that the work in \cite{Petronilho} can be related to the recent method of
computing generating functions via kneading determinants \cite{Alves} applied
to finite and infinite order difference equations. We suggest that this last
method is a good way to tackle the problem of computing explicitly the
solutions of this and other non-autonomous problems in future work.

\subsection{Purpose and overview}

Changing completely the perspective from the above mentioned literature we
approach the problem in the framework of classic linear periodic difference
equations using extensively linear algebra methods.

There exists a simpler method to find the solutions developed by Achille Marie
Gaston Floquet (1847-1920) \cite{Flo1883}. Floquet theory, first published in
1883 for periodic linear differential equations, was extended to difference
equations being a long time classic and familiar in many textbooks
\cite{PFR2004, saberDE, Gol1986, KP2001}. In Floquet theory it is necessary to
find explicitly a monodromy matrix and its eigenvalues, the Floquet multipliers.

We go further since we study the equation (\ref{naFre}) in the general
non-autonomous case for any complex initial conditions and considering
arbitrary complex sequences of parameters $\left\{  a_{n}\right\}  $,
$n=0,1,2,\ldots$, not necessarily periodic. We consider $a_{n}\in\mathbb{%
%TCIMACRO{\U{2102} }%
%BeginExpansion
\mathbb{C}
%EndExpansion
}$ since\ our method works with complex equations and orbits being not
restricted to the original formulation on the natural numbers.

Let us write the non-autonomous difference equation (\ref{naFre}) in matrix
form, namely
\[%
\begin{array}
[c]{ccc}%
\left[
\begin{array}
[c]{c}%
x_{n+2}\\
x_{n+1}%
\end{array}
\right]  = & \left[
\begin{array}
[c]{cc}%
a_{n} & 1\\
1 & 0
\end{array}
\right]  & \left[
\begin{array}
[c]{c}%
x_{n+1}\\
x_{n}%
\end{array}
\right]  ,
\end{array}
\]
for all $n=0,1,2,\ldots$, or equivalently
\begin{equation}
X_{n+1}=A_{n}X_{n}, \label{EqMForm}%
\end{equation}
where
\[
A_{n}=\left[
\begin{array}
[c]{cc}%
a_{n} & 1\\
1 & 0
\end{array}
\right]  ,\text{ }X_{n}=\left[
\begin{array}
[c]{c}%
x_{n+1}\\
x_{n}%
\end{array}
\right]  ,\text{ }X_{n+1}=\left[
\begin{array}
[c]{c}%
x_{n+2}\\
x_{n+1}%
\end{array}
\right]  .
\]
The solution $X_{n}$ of the general non-autonomous case is given by the matrix
product%
\[
X_{n}=A_{n-1}A_{n-2}\ldots A_{1}A_{0}=\left(  \prod_{i=0}^{n-1}A_{n-1-i}%
\right)  X_{0}.
\]
If we call
\begin{equation}
C_{n}=A_{n-1}A_{n-2}\ldots A_{1}A_{0}=\prod_{i=0}^{n-1}A_{n-1-i}
\label{Product}%
\end{equation}
we can write the solution of (\ref{EqMForm}) in the form%
\[
X_{n}=C_{n}X_{0}.
\]
where $X_{0}\in%
%TCIMACRO{\U{2102} }%
%BeginExpansion
\mathbb{C}
%EndExpansion
^{2}$ is a given initial condition. Therefore the core of our work rests on
techniques to obtain the products of the $2\times2$ matrices $A_{n}\in%
%TCIMACRO{\U{2102} }%
%BeginExpansion
\mathbb{C}
%EndExpansion
^{2\times2}$ that are obtained in the main Theorem \ref{Main} of this paper in
section \ref{MR}. It is noticeable that the linear algebra approach is still
complicated regarding that we have only products of very simple $2\times2$ matrices.

When the sequence of numbers $a_{n}$, $n=0,1,2,\dots$ is $k$-periodic with
$a_{k}=a_{0}$ with $k$ minimal, a similar iteration gives $X_{k}=C_{k}X_{0}$.
In the periodic case, the matrix $C_{k}$ is known as the monodromy matrix
\cite{PFR2004, saberDE, Gol1986, KP2001} of the periodic equation
(\ref{EqMForm}). Using the monodromy matrix one can construct the solution of
equation (\ref{EqMForm}) since
\[
X_{mk}=C_{k}^{m}X_{0},
\]
for any $m\geq0$. For a general $n=mk+r$ with $r<k$ we have%
\[
X_{n}=\left(  \prod_{i=0}^{r-1}A_{r-1-i}\right)  C_{k}^{m}X_{0},
\]
and, again, the solution of the periodic case is obtained computing matrix
powers and products of matrices. The key of the periodic problem is to obtain
the eigenvalues of the monodromy matrix, usually called Floquet multipliers in
the field of dynamical systems. Moreover, the asymptotic behavior of the
solutions can be studied via these numbers.

The problem of finding the Floquet multipliers is not easy, since we have to
solve the characteristic equation%
\[
\det\left(  C_{k}-\lambda I\right)  =0
\]
and the entries of $C_{k}$ satisfy the same recurrences of the original
problem. Consequently, in section \ref{MM} we approach this problem computing
directly the eigenvalues of the monodromy matrix, i.e., the Floquet
multipliers of the periodic equation.

In section \ref{Solutions} we study more deeply the solutions of equation
(\ref{naFre}). We determine explicitly conditions for the periodicity or
non-periodicity of the solutions. We remark here that the limits of the
quotients of consecutive iterates exhibit periodicity when the iteration time
$n$ tends to infinity. Something similar to the convergence to the golden
ratio of the quotients of consecutive iterates in the original Fibonacci problem.

Finally, using the techniques previously developed, we present some examples
in section \ref{exampls}. More specifically, we give a complete study of a
$3$-periodic and a $4$-periodic equation where $k=3$ and $k=4$ respectively.

In a nutshell we can say that in this paper the main result rests mainly on
the general technique to compute arbitrary products of some special $2\times2$
matrices. As usual in the field of linear difference equations most part of
the theory can be seen as a reinterpretation of linear algebra results as we
can see in \cite{PFR2004, saberDE, Gol1986, KP2001}. We present explicit
solutions of difference equation (\ref{naFre}) for general sequences of
complex numbers $\left\{  a_{n}\right\}  $, $n=0,1,2,\ldots$, periodic or not,
not obtained before in the literature. We provide abundant examples to help
facilitate clarity.

\section{Main result\label{MR}}

We start this section by introducing some definitions and notations that will
help us greatly in the statement of the results. As usual $\mathbb{%
%TCIMACRO{\U{2115} }%
%BeginExpansion
\mathbb{N}
%EndExpansion
}$ represents the set of non-negative integers and $%
%TCIMACRO{\U{2102} }%
%BeginExpansion
\mathbb{C}
%EndExpansion
\left[  t\right]  $ the set of formal polynomials with complex coefficients in
the indeterminate $t$.

\begin{defn}
Let $\mathcal{A}=\left\{  a_{i}\right\}  _{i\in\mathbb{%
%TCIMACRO{\U{2115} }%
%BeginExpansion
\mathbb{N}
%EndExpansion
}}$ be a non-zero complex sequence of numbers not necessarily periodic. We
define a $2\times2$ \textbf{Fibonacci matrix} $A_{i}\in%
%TCIMACRO{\U{2102} }%
%BeginExpansion
\mathbb{C}
%EndExpansion
^{2\times2}$ as the following%
\[
A_{i}=\left[
\begin{array}
[c]{cc}%
a_{i} & 1\\
1 & 0
\end{array}
\right]  .
\]

\end{defn}

\begin{notation}
Let us use the vector notation $\mathfrak{a}_{i,j}=\left(  a_{i}%
,a_{i+1},\ldots,a_{j-1},a_{j}\right)  ^{T}$, for some $i<j$, $i,j\in\mathbb{%
%TCIMACRO{\U{2115} }%
%BeginExpansion
\mathbb{N}
%EndExpansion
}$. In the sequel, it will be evident that we need to indicate both the first
index and the last index in this notation. When this vector has only one
component the notation will be naturally simplified to $\mathfrak{a}%
_{i}=\left(  \mathfrak{a}_{i}\right)  ^{T}=a_{i}$. Moreover, all the usual
conventions about summations and products will be used.
\end{notation}

\begin{notation}
To simplify the expressions with summations that we have in the following
discussion we define for $n$ and $p$ even the multi summation operation%
\begin{equation}
\widehat{\sum\limits_{i_{0},\ldots,i_{p-1}}^{n,p}}\left(  \cdot\right)
\triangleq\overset{p\text{ terms}}{\overbrace{%
%TCIMACRO{\dsum \limits_{_{i_{0}=0}}^{\frac{n-p}{2}}}%
%BeginExpansion
{\displaystyle\sum\limits_{_{i_{0}=0}}^{\frac{n-p}{2}}}
%EndExpansion%
%TCIMACRO{\dsum \limits_{_{i_{1}=i_{0}}}^{\frac{n-p}{2}}}%
%BeginExpansion
{\displaystyle\sum\limits_{_{i_{1}=i_{0}}}^{\frac{n-p}{2}}}
%EndExpansion%
%TCIMACRO{\dsum \limits_{_{i_{2}=i_{1}+1}}^{\frac{n-p}{2}+1}}%
%BeginExpansion
{\displaystyle\sum\limits_{_{i_{2}=i_{1}+1}}^{\frac{n-p}{2}+1}}
%EndExpansion%
%TCIMACRO{\dsum \limits_{_{i_{3}=i_{2}}}^{\frac{n-p}{2}+1}}%
%BeginExpansion
{\displaystyle\sum\limits_{_{i_{3}=i_{2}}}^{\frac{n-p}{2}+1}}
%EndExpansion
\cdots%
%TCIMACRO{\dsum \limits_{_{i_{p-4}=i_{p-5}+1}}^{\frac{n-2}{2}-1}}%
%BeginExpansion
{\displaystyle\sum\limits_{_{i_{p-4}=i_{p-5}+1}}^{\frac{n-2}{2}-1}}
%EndExpansion%
%TCIMACRO{\dsum \limits_{_{i_{p-3}=i_{p-4}}}^{\frac{n-2}{2}-1}}%
%BeginExpansion
{\displaystyle\sum\limits_{_{i_{p-3}=i_{p-4}}}^{\frac{n-2}{2}-1}}
%EndExpansion%
%TCIMACRO{\dsum \limits_{_{i_{p-2}=i_{p-3}+1}}^{\frac{n-2}{2}}}%
%BeginExpansion
{\displaystyle\sum\limits_{_{i_{p-2}=i_{p-3}+1}}^{\frac{n-2}{2}}}
%EndExpansion%
%TCIMACRO{\dsum \limits_{_{i_{p-1}=i_{p-2}}}^{\frac{n-2}{2}}}%
%BeginExpansion
{\displaystyle\sum\limits_{_{i_{p-1}=i_{p-2}}}^{\frac{n-2}{2}}}
%EndExpansion
}}\left(  \cdot\right)  , \label{MultiSE}%
\end{equation}
when $n$ and $p$ are odd we write%
\begin{equation}
\overset{p\text{ terms}}{\widehat{\sum\limits_{i_{0},\ldots,i_{p-1}}^{n,p}%
}\left(  \cdot\right)  \triangleq\overbrace{%
%TCIMACRO{\dsum \limits_{_{i_{0}=0}}^{\frac{n-p}{2}}}%
%BeginExpansion
{\displaystyle\sum\limits_{_{i_{0}=0}}^{\frac{n-p}{2}}}
%EndExpansion%
%TCIMACRO{\dsum \limits_{_{i_{1}=i_{0}}}^{\frac{n-p}{2}}}%
%BeginExpansion
{\displaystyle\sum\limits_{_{i_{1}=i_{0}}}^{\frac{n-p}{2}}}
%EndExpansion%
%TCIMACRO{\dsum \limits_{_{i_{2}=i_{1}+1}}^{\frac{n-p}{2}+1}}%
%BeginExpansion
{\displaystyle\sum\limits_{_{i_{2}=i_{1}+1}}^{\frac{n-p}{2}+1}}
%EndExpansion%
%TCIMACRO{\dsum \limits_{_{i_{3}=i_{2}}}^{\frac{n-p}{2}+1}}%
%BeginExpansion
{\displaystyle\sum\limits_{_{i_{3}=i_{2}}}^{\frac{n-p}{2}+1}}
%EndExpansion
\cdots%
%TCIMACRO{\dsum \limits_{_{i_{p-3}=i_{p-4}+1}}^{\frac{n-1}{2}-1}}%
%BeginExpansion
{\displaystyle\sum\limits_{_{i_{p-3}=i_{p-4}+1}}^{\frac{n-1}{2}-1}}
%EndExpansion%
%TCIMACRO{\dsum \limits_{_{i_{p-2}=i_{p-3}}}^{\frac{n-1}{2}-1}}%
%BeginExpansion
{\displaystyle\sum\limits_{_{i_{p-2}=i_{p-3}}}^{\frac{n-1}{2}-1}}
%EndExpansion%
%TCIMACRO{\dsum \limits_{_{i_{p-1}=i_{p-2}+1}}^{\frac{n-1}{2}}}%
%BeginExpansion
{\displaystyle\sum\limits_{_{i_{p-1}=i_{p-2}+1}}^{\frac{n-1}{2}}}
%EndExpansion
}}\left(  \cdot\right)  . \label{MultiSO}%
\end{equation}
In both cases the number of sums in each operator is $p$. Please note that the
upper bound on the innermost summation determines the other upper bounds and
not the other way around.
\end{notation}

The introduction of this notation facilitates a great deal the writing of the
proofs since large consecutive sums appear in the entries of the product of
matrices. With this non-standard notation the lower indices and higher indices
appear in consecutive pairs, starting at $0$ for the pair $i_{0}=0$ and
$i_{1}=i_{0}$, having successfully pairs of the type $i_{2m}=i_{2m-1}+1$ and
$i_{2m+1}=i_{2m}$; each pair of upper and lower indices increases by one at
each transition from odd indexed $i_{2m-1\text{ }}$indices to even indexed
$i_{2m}$ indices. Clearly, in the case of odd $n$ and$\ p$ the innermost sum
has only one term instead of a pair of sums. Please check examples
\ref{example1} and \ref{example2}, or the last sections of this paper, to see
practical evaluations.

\begin{defn}
Let $n\in\mathbb{%
%TCIMACRO{\U{2115} }%
%BeginExpansion
\mathbb{N}
%EndExpansion
}$ be even. We write $\chi_{n,p}\left(  \mathfrak{a}_{0,n-1}\right)  t^{p}\in%
%TCIMACRO{\U{2102} }%
%BeginExpansion
\mathbb{C}
%EndExpansion
\left[  t\right]  $, with even degree $p$ such that $0<p\leq n$, as the
following formal monomial in the indeterminate $t$
\begin{equation}
\chi_{n,p}\left(  \mathfrak{a}_{0,n-1}\right)  t^{p}\triangleq t^{p}%
\widehat{\sum\limits_{i_{0},\ldots,i_{p-1}}^{n,p}}%
%TCIMACRO{\dprod _{l=0}^{\frac{p}{2}-1}}%
%BeginExpansion
{\displaystyle\prod_{l=0}^{\frac{p}{2}-1}}
%EndExpansion
a_{2i_{2l}}%
%TCIMACRO{\dprod _{l=0}^{\frac{p}{2}-1}}%
%BeginExpansion
{\displaystyle\prod_{l=0}^{\frac{p}{2}-1}}
%EndExpansion
a_{1+2i_{2l+1}}. \label{CHI}%
\end{equation}
Finally, when we do not restrict the indices in $\mathfrak{a}$ to start at $0$
but at $L\in%
%TCIMACRO{\U{2115} }%
%BeginExpansion
\mathbb{N}
%EndExpansion
$, we have $\chi_{n,p}\left(  \mathfrak{a}_{L,L+n-1}\right)  $ such that
\[
\chi_{n,p}\left(  \mathfrak{a}_{L,L+n-1}\right)  t^{p}\triangleq
t^{p}\widehat{\sum\limits_{i_{0},\ldots,i_{p-1}}^{n,p}}%
%TCIMACRO{\dprod _{l=0}^{\frac{p}{2}-1}}%
%BeginExpansion
{\displaystyle\prod_{l=0}^{\frac{p}{2}-1}}
%EndExpansion
a_{L+2i_{2l}}%
%TCIMACRO{\dprod _{l=0}^{\frac{p}{2}-1}}%
%BeginExpansion
{\displaystyle\prod_{l=0}^{\frac{p}{2}-1}}
%EndExpansion
a_{L+1+2i_{2l+1}}.
\]
When $p=0$, we consider by definition that $\chi_{n,0}\left(  v\right)
\triangleq1$, for every even $n$ and any vector $v$. Moreover, $\chi
_{0,0}\left(  \varnothing\right)  =\chi_{0,0}=1$.
\end{defn}

The apparently complicated expression (\ref{CHI})\ is nothing but the sum of
all possible products $a_{j_{1}}\ldots a_{j_{p}}$ with $p$ factors, such that
the first factor $a_{j_{1}}$ has even index $j_{1}$, the second factor has odd
index $j_{2}$ and so on, ending each product with an odd index factor
$a_{j_{p}}$ where $j_{p}$ is odd.

\begin{exmp}
\label{example1}One can see that
\begin{align*}
\chi_{6,4}\left(  \mathfrak{a}_{0,5}\right)   &  =\widehat{\sum\limits_{i_{0}%
,i_{1},i_{2},i_{3}}^{6,4}}\prod\limits_{l=0}^{1}a_{2i_{2l}}\prod
\limits_{l=0}^{1}a_{1+2i_{2l+1}}\\
& \\
&  =\sum_{i_{0}=0}^{1}\sum_{i_{1}=i_{0}}^{1}\sum_{i_{2}=i_{1}+1}^{2}%
\sum_{i_{3}=i_{2}}^{2}\prod\limits_{l=0}^{1}a_{2i_{2l}}\prod\limits_{l=0}%
^{1}a_{1+2i_{2l+1}}\\
& \\
&  =a_{0}a_{1}a_{2}a_{3}+a_{0}a_{1}a_{2}a_{5}+a_{0}a_{1}a_{4}a_{5}+a_{0}%
a_{3}a_{4}a_{5}+a_{2}a_{3}a_{4}a_{5}.
\end{align*}
Obviously, when we shift the indices in $\mathfrak{a}$ by one unit, we have%
\[
\chi_{6,4}\left(  \mathfrak{a}_{1,6}\right)  =a_{1}a_{2}a_{3}a_{4}+a_{1}%
a_{2}a_{3}a_{6}+a_{1}a_{2}a_{5}a_{6}+a_{1}a_{4}a_{5}a_{6}+a_{3}a_{4}a_{5}%
a_{6}.
\]
From the practical point of view, one can see that the construction of the
number $\chi_{6,4}\left(  \mathfrak{a}_{0,5}\right)  $ is equivalent to the
combinatorial problem of finding all the possible configurations obtained from
$\left(  0,1,2,3,4,5\right)  $ when we cut a pair of consecutive numbers (not
considering $a_{5}$ and $a_{0}$ consecutive) in the sequence of elements of
the set $\mathcal{A}$. In this case, we have the configurations
\[
\left\{  \left(  0,1,2,3\right)  ,\left(  0,1,2,5\right)  ,\left(
0,1,4,5\right)  ,\left(  0,3,4,5\right)  ,\left(  2,3,4,5\right)  \right\}  .
\]

\end{exmp}

One can formulate this problem, for example, as the possible configurations of
people remaining when a pair of persons seating together leave a counter or
the possible configurations of a row of distinguishable balls when two
adjacent balls are taken.

\begin{defn}
Consider $n$ odd. Similarly to the even case, we define the $p$-degree formal
monomial $\chi_{n,p}\left(  \mathfrak{a}_{L,L+n-1}\right)  t^{p}\in%
%TCIMACRO{\U{2102} }%
%BeginExpansion
\mathbb{C}
%EndExpansion
\left[  t\right]  $, with $1\leq p\leq n$, as
\[
\chi_{n,p}\left(  \mathfrak{a}_{L,L+n-1}\right)  t^{p}\triangleq
t^{p}\widehat{\sum\limits_{i_{0},\ldots,i_{p-1}}^{n,p}}%
%TCIMACRO{\dprod _{l=0}^{\frac{p-1}{2}}}%
%BeginExpansion
{\displaystyle\prod_{l=0}^{\frac{p-1}{2}}}
%EndExpansion
a_{L+2i_{2l}}%
%TCIMACRO{\dprod _{l=0}^{\frac{p-1}{2}-1}}%
%BeginExpansion
{\displaystyle\prod_{l=0}^{\frac{p-1}{2}-1}}
%EndExpansion
a_{L+1+2i_{2l+1}}.
\]

\end{defn}

We remark again that the apparently complicated expression of $\chi
_{n,p}\left(  \mathfrak{a}_{0,n-1}\right)  $ is the sum of all possible
products $a_{j_{1}}\ldots a_{j_{p}}$ with $p$ factors, such that the first
factor has even index, the second factor has odd index and so on, ending on
$a_{j_{p}}$ where $j_{p}$ is an even index.

\begin{exmp}
\label{example2}For instance, one can see that
\begin{align*}
\chi_{5,3}\left(  \mathfrak{a}_{0,4}\right)   &  =\widehat{\sum\limits_{i_{0}%
,i_{1},i_{2}}^{5,3}}\prod\limits_{l=0}^{1}a_{2i_{2l}}\prod\limits_{l=0}%
^{0}a_{1+2i_{2l+1}}\\
& \\
&  =\sum_{i_{0}=0}^{1}\sum_{i_{1}=i_{0}}^{1}\sum_{i_{2}=i_{1}+1}^{2}%
\prod\limits_{l=0}^{1}a_{2i_{2l}}\prod\limits_{l=0}^{0}a_{1+2i_{2l+1}}\\
& \\
&  =a_{0}a_{1}a_{2}+a_{0}a_{1}a_{4}+a_{0}a_{3}a_{4}+a_{2}a_{3}a_{4}.
\end{align*}
Again, this is the same combinatorial problem of finding all the possible
configurations obtained from $\left(  0,1,2,3,4\right)  $ when we cut a pair
of consecutive numbers. In this case the possible configurations are%
\[
\left\{  \left(  0,1,2\right)  ,\left(  0,1,4\right)  ,\left(  0,3,4\right)
,\left(  2,3,4\right)  \right\}  .
\]

\end{exmp}

Other easy cases are $\chi_{5,1}\left(  \mathfrak{a}_{0,4}\right)
=a_{0}+a_{2}+a_{4}$ or $\chi_{5,1}\left(  \mathfrak{a}_{1,5}\right)
=a_{1}+a_{3}+a_{5}$.

\begin{defn}
For the natural numbers $n\geq0$ and $i\in%
%TCIMACRO{\U{2115} }%
%BeginExpansion
\mathbb{N}
%EndExpansion
$, we define the formal polynomial $\Omega\left(  t\right)  \in%
%TCIMACRO{\U{2102} }%
%BeginExpansion
\mathbb{C}
%EndExpansion
\left[  t\right]  $ by%
\begin{equation}
\Omega_{n}\left(  t,\mathfrak{a}_{L,L+n-1}\right)  \triangleq\left\{
\begin{array}
[c]{l}%
%TCIMACRO{\dsum \limits_{j=0}^{\frac{n}{2}}}%
%BeginExpansion
{\displaystyle\sum\limits_{j=0}^{\frac{n}{2}}}
%EndExpansion
\chi_{n,2j}\left(  \mathfrak{a}_{L,L+n-1}\right)  t^{2j}\text{, }n\text{
even,}\\
\\%
%TCIMACRO{\dsum \limits_{j=0}^{\frac{n-1}{2}}}%
%BeginExpansion
{\displaystyle\sum\limits_{j=0}^{\frac{n-1}{2}}}
%EndExpansion
\chi_{n,2j+1}\left(  \mathfrak{a}_{L,L+n-1}\right)  t^{2j+1}\text{, }n\text{
odd.}%
\end{array}
\right.  \label{FP_OMEGA}%
\end{equation}

\end{defn}

The polynomial $\Omega_{n}\left(  t,\mathfrak{a}_{0,n-1}\right)  $ has degree
$n$ and when $n=0$ we write $\Omega_{0}\left(  t\right)  =\chi_{0,0}t^{0}=1$.
The construction of $\Omega_{n}\left(  t,\mathfrak{a}_{0,n-1}\right)  $
corresponds to the combinatorial problem of finding all the possible
configurations when we cut all the possible pairs, from zero pairs to
$\left\lfloor \frac{n}{2}\right\rfloor $ pairs, of consecutive numbers from a
row of $n$ integers.

\begin{exmp}
As an example, let us illustrate how to compute $\Omega_{2}\left(
t,\mathfrak{a}_{0,1}\right)  $. Since
\[
\Omega_{2}\left(  t,\mathfrak{a}_{0,1}\right)  =\chi_{2,0}\left(
\mathfrak{a}_{0,1}\right)  +\chi_{2,2}\left(  \mathfrak{a}_{0,1}\right)
t^{2},
\]
it follows that
\[
\chi_{2,2}\left(  \mathfrak{a}_{0,1}\right)  =\widehat{\sum\limits_{i_{0}%
,i_{1}}^{2,2}}%
%TCIMACRO{\dprod _{l=0}^{0}}%
%BeginExpansion
{\displaystyle\prod_{l=0}^{0}}
%EndExpansion
a_{2i_{2l}}%
%TCIMACRO{\dprod _{l=0}^{0}}%
%BeginExpansion
{\displaystyle\prod_{l=0}^{0}}
%EndExpansion
a_{1+2i_{2l+1}}=a_{0}a_{1},
\]
and
\[
\chi_{2,0}\left(  \mathfrak{a}_{0,1}\right)  =1.
\]
Hence, $\Omega_{2}\left(  t,\mathfrak{a}_{0,1}\right)  =1+a_{0}a_{1}t^{2}$.
Notice that $\Omega_{1}\left(  t,\mathfrak{a}_{0}\right)  =a_{0}t$ and
$\Omega_{0}\left(  t\right)  =1$. Finally we note that in the case of $i=1$
the indices in $\chi_{2,2}\left(  \mathfrak{a}_{0,1}\right)  $ are increased
by $1$, so $\Omega_{2}\left(  t,\mathfrak{a}_{1,2}\right)  =1+a_{1}a_{2}t^{2}$
\end{exmp}

In the next result we show how to compute the product matrix $C_{n}$ defined
in (\ref{Product}).

\begin{thm}
\label{Main}Consider the formal polynomials $\Omega_{n}\left(  t,\mathfrak{a}%
_{0,n-1}\right)  $, $\Omega_{n-1}\left(  t,\mathfrak{a}_{1,n-1}\right)  $,
$\Omega_{n-1}\left(  t,\mathfrak{a}_{0,n-2}\right)  $, $\Omega_{n-2}\left(
t,\mathfrak{a}_{1,n-2}\right)  \in%
%TCIMACRO{\U{2102} }%
%BeginExpansion
\mathbb{C}
%EndExpansion
\left[  t\right]  $ defined in (\ref{FP_OMEGA}). For any $n\geq2$ and $t=1$,
the product matrix of Fibonacci matrices defined in (\ref{Product}) is given
by
\[
C_{n}=\prod_{i=0}^{n-1}\left[
\begin{array}
[c]{cc}%
a_{n-1-i} & 1\\
1 & 0
\end{array}
\right]  =\left[
\begin{array}
[c]{cc}%
\Omega_{n}\left(  1,\mathfrak{a}_{0,n-1}\right)  & \Omega_{n-1}\left(
1,\mathfrak{a}_{1,n-1}\right) \\
\Omega_{n-1}\left(  1,\mathfrak{a}_{0,n-2}\right)  & \Omega_{n-2}\left(
1,\mathfrak{a}_{1,n-2}\right)
\end{array}
\right]  .
\]

\end{thm}

\begin{pf}
The proof is by induction in $n$. When $n=2$ we have
\[
C_{2}=\left[
\begin{array}
[c]{cc}%
\Omega_{2}\left(  1,\mathfrak{a}_{0,1}\right)  & \Omega_{1}\left(
1,\mathfrak{a}_{1}\right) \\
\Omega_{1}\left(  1,\mathfrak{a}_{0}\right)  & \Omega_{0}\left(  t\right)
\end{array}
\right]  =\left[
\begin{array}
[c]{cc}%
1+a_{0}a_{1} & a_{1}\\
a_{0} & 1
\end{array}
\right]  ,
\]
which is $%
%TCIMACRO{\dprod \limits_{j=0}^{1}}%
%BeginExpansion
{\displaystyle\prod\limits_{j=0}^{1}}
%EndExpansion
A_{1-j}$. We notice also that%
\[
C_{3}=\left[
\begin{array}
[c]{cc}%
\Omega_{3}\left(  1,\mathfrak{a}_{0,2}\right)  & \Omega_{2}\left(
1,\mathfrak{a}_{1,2}\right) \\
\Omega_{2}\left(  1,\mathfrak{a}_{0,1}\right)  & \Omega_{1}\left(
1,a_{1}\right)
\end{array}
\right]  =\left[
\begin{array}
[c]{cc}%
a_{0}+a_{2}+a_{0}a_{1}a_{2} & 1+a_{1}a_{2}\\
1+a_{0}a_{1} & a_{1}%
\end{array}
\right]  .
\]
It is an easy computation to show that this matrix is $%
%TCIMACRO{\dprod \limits_{j=0}^{2}}%
%BeginExpansion
{\displaystyle\prod\limits_{j=0}^{2}}
%EndExpansion
A_{2-j}$.

Now we have the hypothesis%
\[
C_{n}=\left[
\begin{array}
[c]{cc}%
\Omega_{n}\left(  1,\mathfrak{a}_{0,n-1}\right)  & \Omega_{n-1}\left(
1,\mathfrak{a}_{1,n-1}\right) \\
\Omega_{n-1}\left(  1,\mathfrak{a}_{0,n-2}\right)  & \Omega_{n-2}\left(
1,\mathfrak{a}_{1,n-2}\right)
\end{array}
\right]  .
\]
One needs to show the induction step for the matrix $C_{n+1}$,
\begin{align*}
C_{n+1}  &  =\left[
\begin{array}
[c]{cc}%
c_{n+1}\left(  1,1\right)  & c_{n+1}\left(  1,2\right) \\
c_{n+1}\left(  2,1\right)  & c_{n+1}\left(  2,2\right)
\end{array}
\right] \\
&  =\left[
\begin{array}
[c]{cc}%
a_{n} & 1\\
1 & 0
\end{array}
\right]  C_{n}%
\end{align*}%
\[
=\left[
\begin{array}
[c]{cc}%
^{a_{n}\Omega_{n}\left(  1,\mathfrak{a}_{0,n-1}\right)  +\Omega_{n-1}\left(
1,\mathfrak{a}_{0,n-2}\right)  } & ^{a_{n}\Omega_{n-1}\left(  1,\mathfrak{a}%
_{1,n-1}\right)  +\Omega_{n-2}\left(  1,\mathfrak{a}_{1,n-2}\right)  }\\
^{\Omega_{n}\left(  1,\mathfrak{a}_{0,n-1}\right)  } & ^{\Omega_{n-1}\left(
1,\mathfrak{a}_{1,n-1}\right)  }%
\end{array}
\right]  .
\]
Automatically, the entries $c_{n+1}\left(  2,1\right)  $ and $c_{n+1}\left(
2,2\right)  $ are done. Now, we have to prove that%
\begin{align}
\Omega_{n+1}\left(  1,\mathfrak{a}_{0,n}\right)   &  =c_{n+1}\left(
1,1\right) \label{caso n+1}\\
&  =a_{n}\Omega_{n}\left(  1,\mathfrak{a}_{0,n-1}\right)  +\Omega_{n-1}\left(
1,\mathfrak{a}_{0,n-2}\right)  ,\nonumber
\end{align}
for any $n$.

We note that%
\begin{align*}
\Omega_{n}\left(  1,\mathfrak{a}_{0,n}\right)   &  =c_{n+1}\left(  1,2\right)
\\
&  =a_{n}\Omega_{n-1}\left(  1,\mathfrak{a}_{1,n-1}\right)  +\Omega
_{n-2}\left(  1,\mathfrak{a}_{1,n-2}\right)  ,
\end{align*}
holds immediately if the condition (\ref{caso n+1})\ is true.

We make use again of the formal indeterminate $t$ to tackle the polynomial%
\[
a_{n}t\Omega_{n}\left(  t,\mathfrak{a}_{0,n-1}\right)  +\Omega_{n-1}\left(
t,\mathfrak{a}_{0,n-2}\right)  ,
\]
computed for each degree of $t$. At the end of the proof, we will make $t=1$
to obtain the desired equality.

Supposing first that $n$ is even we note that%
\begin{align*}
a_{n}t\Omega_{n}\left(  t,\mathfrak{a}_{0,n-1}\right)  +\Omega_{n-1}\left(
t,\mathfrak{a}_{0,n-2}\right)   &  =\\
&
\end{align*}%
\begin{equation}
a_{n}t%
%TCIMACRO{\dsum \limits_{j=0}^{\frac{n}{2}}}%
%BeginExpansion
{\displaystyle\sum\limits_{j=0}^{\frac{n}{2}}}
%EndExpansion
\chi_{n,2j}\left(  \mathfrak{a}_{0,n-1}\right)  t^{2j}+%
%TCIMACRO{\dsum \limits_{j=0}^{\frac{n}{2}}}%
%BeginExpansion
{\displaystyle\sum\limits_{j=0}^{\frac{n}{2}}}
%EndExpansion
\chi_{n-1,2j+1}\left(  \mathfrak{a}_{0,n-2}\right)  t^{2j+1}. \label{even}%
\end{equation}
The term with the highest degree is
\begin{align*}
a_{n}t\chi_{n,n}\left(  \mathfrak{a}_{0,n-1}\right)  t^{n}  &  =t^{n+1}a_{n}%
%TCIMACRO{\dprod _{l=0}^{n-1}}%
%BeginExpansion
{\displaystyle\prod_{l=0}^{n-1}}
%EndExpansion
a_{l}\\
& \\
&  =t^{n+1}%
%TCIMACRO{\dprod _{l=0}^{n}}%
%BeginExpansion
{\displaystyle\prod_{l=0}^{n}}
%EndExpansion
a_{l}=\chi_{n+1,n+1}\left(  \mathfrak{a}_{0,n}\right)  t^{n+1}.
\end{align*}
On the other hand, the term with lowest degree is given by
\begin{align*}
a_{n}t\chi_{n,0}\left(  \mathfrak{a}_{0,n-1}\right)  +\chi_{n-1,1}\left(
\mathfrak{a}_{0,n-2}\right)  t  &  =a_{n}t+%
%TCIMACRO{\dsum \limits_{i=0}^{n-1}}%
%BeginExpansion
{\displaystyle\sum\limits_{i=0}^{n-1}}
%EndExpansion
a_{i}t=%
%TCIMACRO{\dsum \limits_{i=0}^{n}}%
%BeginExpansion
{\displaystyle\sum\limits_{i=0}^{n}}
%EndExpansion
a_{i}t\\
& \\
&  =\chi_{n+1,1}\left(  \mathfrak{a}_{0,n}\right)  t.\medskip
\end{align*}
Please note that $n$ is even and $p$ must also be even. We consider the sums
in the same degree $p$ in (\ref{even})
\[
a_{n}t\chi_{n,p}\left(  \mathfrak{a}_{0,n-1}\right)  t^{p}+\chi_{n-1,p+1}%
\left(  \mathfrak{a}_{0,n-2}\right)  t^{p+1}.
\]
Consequently%
\[
a_{n}\chi_{n,p}\left(  \mathfrak{a}_{0,n-1}\right)  +\chi_{n-1,p+1}\left(
\mathfrak{a}_{0,n-2}\right)  =\medskip
\]%
\[
\widehat{\sum\limits_{i_{0},\ldots,i_{p-1}}^{n,p}}a_{n}%
%TCIMACRO{\dprod _{l=0}^{\frac{p}{2}-1}}%
%BeginExpansion
{\displaystyle\prod_{l=0}^{\frac{p}{2}-1}}
%EndExpansion
a_{2i_{2l}}%
%TCIMACRO{\dprod _{l=0}^{\frac{p}{2}-1}}%
%BeginExpansion
{\displaystyle\prod_{l=0}^{\frac{p}{2}-1}}
%EndExpansion
a_{1+2i_{2l+1}}+\widehat{\sum\limits_{i_{0},\ldots,i_{p}}^{n-1,p+1}}%
%TCIMACRO{\dprod _{l=0}^{\frac{p}{2}}}%
%BeginExpansion
{\displaystyle\prod_{l=0}^{\frac{p}{2}}}
%EndExpansion
a_{2i_{2l}}%
%TCIMACRO{\dprod _{l=0}^{\frac{p}{2}-1}}%
%BeginExpansion
{\displaystyle\prod_{l=0}^{\frac{p}{2}-1}}
%EndExpansion
a_{1+2i_{2l+1}}.\medskip
\]
Consider now%
\[
\chi_{n+1,p+1}\left(  \mathfrak{a}_{0,n}\right)  =\widehat{\sum\limits_{i_{0}%
,\ldots,i_{p}}^{n-1,p+1}}%
%TCIMACRO{\dprod _{l=0}^{\frac{p}{2}}}%
%BeginExpansion
{\displaystyle\prod_{l=0}^{\frac{p}{2}}}
%EndExpansion
a_{2i_{2l}}%
%TCIMACRO{\dprod _{l=0}^{\frac{p}{2}-1}}%
%BeginExpansion
{\displaystyle\prod_{l=0}^{\frac{p}{2}-1}}
%EndExpansion
a_{1+2i_{2l+1}}.\medskip
\]
We split this last member in two terms, one with factor $a_{n}$ and the other
without $a_{n}$. Hence,
\[
\chi_{n+1,p+1}\left(  \mathfrak{a}_{0,n}\right)  =S_{1}+S_{2},
\]
where
\[
S_{1}=\widehat{\sum\limits_{i_{0},\ldots,i_{p-1}}^{n,p}}%
%TCIMACRO{\dsum \limits_{i_{p}=i_{p-1}+1}^{\frac{n}{2}-1}}%
%BeginExpansion
{\displaystyle\sum\limits_{i_{p}=i_{p-1}+1}^{\frac{n}{2}-1}}
%EndExpansion%
%TCIMACRO{\dprod _{l=0}^{\frac{p}{2}}}%
%BeginExpansion
{\displaystyle\prod_{l=0}^{\frac{p}{2}}}
%EndExpansion
a_{2i_{2l}}%
%TCIMACRO{\dprod _{l=0}^{\frac{p}{2}-1}}%
%BeginExpansion
{\displaystyle\prod_{l=0}^{\frac{p}{2}-1}}
%EndExpansion
a_{1+2i_{2l+1}}%
\]
and
\begin{align*}
S_{2}  &  =\widehat{\sum\limits_{i_{0},\ldots,i_{p-1}}^{n,p}}a_{n}%
%TCIMACRO{\dprod _{l=0}^{\frac{p}{2}-1}}%
%BeginExpansion
{\displaystyle\prod_{l=0}^{\frac{p}{2}-1}}
%EndExpansion
a_{2i_{2l}}%
%TCIMACRO{\dprod _{l=0}^{\frac{p}{2}-1}}%
%BeginExpansion
{\displaystyle\prod_{l=0}^{\frac{p}{2}-1}}
%EndExpansion
a_{1+2i_{2l+1}}\\
& \\
&  =a_{n}\chi_{n,p}\left(  \mathfrak{a}_{0,n-1}\right)  .
\end{align*}
On the other hand%
\begin{equation}
S_{1}=\widehat{\sum\limits_{i_{0},\ldots,i_{p-1}}^{n,p}}%
%TCIMACRO{\dsum \limits_{i_{p}=i_{p-1}+1}^{\frac{n}{2}-1}}%
%BeginExpansion
{\displaystyle\sum\limits_{i_{p}=i_{p-1}+1}^{\frac{n}{2}-1}}
%EndExpansion%
%TCIMACRO{\dprod _{l=0}^{\frac{p}{2}}}%
%BeginExpansion
{\displaystyle\prod_{l=0}^{\frac{p}{2}}}
%EndExpansion
a_{2i_{2l}}%
%TCIMACRO{\dprod _{l=0}^{\frac{p}{2}-1}}%
%BeginExpansion
{\displaystyle\prod_{l=0}^{\frac{p}{2}-1}}
%EndExpansion
a_{1+2i_{2l+1}}.\medskip\label{S1}%
\end{equation}
Now $p+1$ is odd. The innermost summation upper bound forces all the other
upper bounds in all the summations to be decreasing one unit by steps of two
summations. The nominal upper bounds do not matter if bigger than the forced
upper bounds since the summations cannot be done due to the lack of possible
summands when the upper bounds exceed the forced ones. The nominal upper
bounds of all summations $\widehat{\sum\limits_{i_{0},\ldots,i_{p-1}}^{n,p}}%
%TCIMACRO{\dsum \limits_{i_{p}=i_{p-1}+1}^{\frac{n}{2}-1}}%
%BeginExpansion
{\displaystyle\sum\limits_{i_{p}=i_{p-1}+1}^{\frac{n}{2}-1}}
%EndExpansion
$ in (\ref{S1}), from inner summation to outer summation, are%
\[
\overset{p+1}{\overbrace{\frac{n}{2}-1,\frac{n}{2}-1,\frac{n}{2}%
-1,\ldots,\frac{n-p}{2},\frac{n-p}{2}}}%
\]
therefore the actual upper bounds that really matter for the computation, from
innermost summation to the outermost summation, are%
\[
\overset{p+1}{\overbrace{\frac{n}{2}-1,\frac{n}{2}-2,\frac{n}{2}%
-2,\ldots,\frac{n-p-2}{2},\frac{n-p-2}{2}}}%
\]
corresponding to odd values $n-1$ and $p+1$ in expression (\ref{MultiSO}).
This leads to the conclusion that%
\begin{align*}
S_{1}  &  =\widehat{\sum\limits_{i_{0},\ldots,i_{p}}^{n-1,p+1}}%
%TCIMACRO{\dprod _{l=0}^{\frac{p}{2}}}%
%BeginExpansion
{\displaystyle\prod_{l=0}^{\frac{p}{2}}}
%EndExpansion
a_{2i_{2l}}%
%TCIMACRO{\dprod _{l=0}^{\frac{p}{2}-1}}%
%BeginExpansion
{\displaystyle\prod_{l=0}^{\frac{p}{2}-1}}
%EndExpansion
a_{1+2i_{2l+1}}\\
& \\
&  =\chi_{n-1,p+1}\left(  \mathfrak{a}_{0,n-2}\right)  \text{.}%
\end{align*}

Adding all the possible values of $\chi_{n+1,2j+1}\left(  \mathfrak{a}%
_{0,n}\right)  t^{2j+1}$ and putting $t=1$ we get%
\[
\Omega_{n+1}\left(  1,\mathfrak{a}_{0,n}\right)  =a_{n}\Omega_{n}\left(
1,\mathfrak{a}_{0,n-1}\right)  +\Omega_{n-1}\left(  1,\mathfrak{a}%
_{0,n-2}\right)  .
\]
The reasonings for $n$ odd and $p$ odd are exactly the same. For sake of
completeness we present here the corresponding induction step for odd $n$.

Let $n$ be odd. Then, one can verify that%
\[
a_{n}t\Omega_{n}\left(  t,\mathfrak{a}_{0,n-1}\right)  +\Omega_{n-1}\left(
t,\mathfrak{a}_{0,n-2}\right)  =\vspace{0.5cm}%
\]%
\begin{equation}
a_{n}t%
%TCIMACRO{\dsum \limits_{j=0}^{\frac{n-1}{2}}}%
%BeginExpansion
{\displaystyle\sum\limits_{j=0}^{\frac{n-1}{2}}}
%EndExpansion
\chi_{n,2j+1}\left(  \mathfrak{a}_{0,n-1}\right)  t^{2j+1}+%
%TCIMACRO{\dsum \limits_{j=0}^{\frac{n-1}{2}}}%
%BeginExpansion
{\displaystyle\sum\limits_{j=0}^{\frac{n-1}{2}}}
%EndExpansion
\chi_{n-1,2j}\left(  \mathfrak{a}_{0,n-2}\right)  t^{2j}. \label{odd}%
\end{equation}
The term with the highest degree is
\begin{align*}
a_{n}t\chi_{n,n}\left(  \mathfrak{a}_{0,n-1}\right)  t^{n}  &  =t^{n+1}a_{n}%
%TCIMACRO{\dprod _{l=0}^{n-1}}%
%BeginExpansion
{\displaystyle\prod_{l=0}^{n-1}}
%EndExpansion
a_{l}=t^{n+1}%
%TCIMACRO{\dprod _{l=0}^{n}}%
%BeginExpansion
{\displaystyle\prod_{l=0}^{n}}
%EndExpansion
a_{l}\\
& \\
&  =\chi_{n+1,n+1}\left(  \mathfrak{a}_{0,n}\right)  t^{n+1},
\end{align*}
while the term with lowest degree is
\[
\chi_{n-1,0}\left(  \mathfrak{a}_{0,n-2}\right)  t^{0}=1=\chi_{n+1,0}\left(
\mathfrak{a}_{0,n}\right)  .
\]
Let us consider the sums with the same degree $p$ in (\ref{odd})%
\[
a_{n}t\chi_{n,p}\left(  \mathfrak{a}_{0,n-1}\right)  t^{p}+\chi_{n-1,p+1}%
\left(  \mathfrak{a}_{0,n-2}\right)  t^{p+1}.
\]
Hence, we have to study the coefficients of $t^{p+1}$%
\[
a_{n}\chi_{n,p}\left(  \mathfrak{a}_{0,n-1}\right)  +\chi_{n-1,p+1}\left(
\mathfrak{a}_{0,n-2}\right)  =
\]%
\[
\widehat{\sum\limits_{i_{0},\ldots,i_{p-1}}^{n,p}}a_{n}%
%TCIMACRO{\dprod _{l=0}^{\frac{p-1}{2}}}%
%BeginExpansion
{\displaystyle\prod_{l=0}^{\frac{p-1}{2}}}
%EndExpansion
a_{2i_{2l}}%
%TCIMACRO{\dprod _{l=0}^{\frac{p-1}{2}-1}}%
%BeginExpansion
{\displaystyle\prod_{l=0}^{\frac{p-1}{2}-1}}
%EndExpansion
a_{1+2i_{2l+1}}+\widehat{\sum\limits_{i_{0},\ldots,i_{p}}^{n-1,p+1}}%
%TCIMACRO{\dprod _{l=0}^{\frac{p-1}{2}}}%
%BeginExpansion
{\displaystyle\prod_{l=0}^{\frac{p-1}{2}}}
%EndExpansion
a_{2i_{2l}}%
%TCIMACRO{\dprod _{l=0}^{\frac{p-1}{2}}}%
%BeginExpansion
{\displaystyle\prod_{l=0}^{\frac{p-1}{2}}}
%EndExpansion
a_{1+2i_{2l+1}}\text{.}%
\]
Consider now%
\[
\chi_{n+1,p+1}\left(  \mathfrak{a}_{0,n}\right)  =\widehat{\sum\limits_{i_{0}%
,\ldots,i_{p}}^{n+1,p+1}}%
%TCIMACRO{\dprod _{l=0}^{\frac{p-1}{2}}}%
%BeginExpansion
{\displaystyle\prod_{l=0}^{\frac{p-1}{2}}}
%EndExpansion
a_{2\ i_{2l}}%
%TCIMACRO{\dprod _{l=0}^{\frac{p-1}{2}}}%
%BeginExpansion
{\displaystyle\prod_{l=0}^{\frac{p-1}{2}}}
%EndExpansion
a_{1+2\ i_{2l+1}}.
\]
Splitting this last member in two terms, one with factor $a_{n}$ and the other
without this factor we have
\[
\chi_{n+1,p+1}\left(  \mathfrak{a}_{0,n}\right)  =\widetilde{S}_{1}%
+\widetilde{S}_{2},
\]
where%
\[
\widetilde{S}_{1}=\widehat{\sum\limits_{i_{0},\ldots,i_{p-1}}^{n,p}}%
%TCIMACRO{\dsum \limits_{i_{p}=i_{p-1}}^{\frac{n-1}{2}-1}}%
%BeginExpansion
{\displaystyle\sum\limits_{i_{p}=i_{p-1}}^{\frac{n-1}{2}-1}}
%EndExpansion%
%TCIMACRO{\dprod _{l=0}^{\frac{p-1}{2}}}%
%BeginExpansion
{\displaystyle\prod_{l=0}^{\frac{p-1}{2}}}
%EndExpansion
a_{2\ i_{2l}}%
%TCIMACRO{\dprod _{l=0}^{\frac{p-1}{2}}}%
%BeginExpansion
{\displaystyle\prod_{l=0}^{\frac{p-1}{2}}}
%EndExpansion
a_{1+2\ i_{2l+1}}%
\]
and%
\[
\widetilde{S}_{2}=\widehat{\sum\limits_{i_{0},\ldots,i_{p-1}}^{n,p}}a_{n}%
%TCIMACRO{\dprod _{l=0}^{\frac{p-1}{2}}}%
%BeginExpansion
{\displaystyle\prod_{l=0}^{\frac{p-1}{2}}}
%EndExpansion
a_{2\ i_{2l}}%
%TCIMACRO{\dprod _{l=0}^{\frac{p-1}{2}-1}}%
%BeginExpansion
{\displaystyle\prod_{l=0}^{\frac{p-1}{2}-1}}
%EndExpansion
a_{1+2\ i_{2l+1}}.
\]
One can verify that $\widetilde{S}_{2}=a_{n}\chi_{n,p}\left(  \mathfrak{a}%
_{0,n-1}\right)  $ and $\widetilde{S}_{1}$ can be written as%
\[
\widetilde{S}_{1}=\widehat{\sum\limits_{i_{0},\ldots,i_{p-1}}^{n,p}}%
%TCIMACRO{\dsum \limits_{i_{p}=i_{p-1}}^{\frac{n-1-2}{2}}}%
%BeginExpansion
{\displaystyle\sum\limits_{i_{p}=i_{p-1}}^{\frac{n-1-2}{2}}}
%EndExpansion%
%TCIMACRO{\dprod \limits_{l=0}^{\frac{p-1}{2}}}%
%BeginExpansion
{\displaystyle\prod\limits_{l=0}^{\frac{p-1}{2}}}
%EndExpansion
a_{2i_{2l}}%
%TCIMACRO{\dprod \limits_{l=0}^{\frac{p-1}{2}}}%
%BeginExpansion
{\displaystyle\prod\limits_{l=0}^{\frac{p-1}{2}}}
%EndExpansion
a_{1+2i_{2l+1}},
\]
a similar reasoning to the one used to simplify (\ref{S1}) gives%
\begin{align*}
\widetilde{S}_{1}  &  =\widehat{\sum\limits_{i_{0},\ldots,i_{p}}^{n-1,p+1}}%
%TCIMACRO{\dprod \limits_{l=0}^{\frac{p-1}{2}}}%
%BeginExpansion
{\displaystyle\prod\limits_{l=0}^{\frac{p-1}{2}}}
%EndExpansion
a_{2i_{2l}}%
%TCIMACRO{\dprod \limits_{l=0}^{\frac{p-1}{2}}}%
%BeginExpansion
{\displaystyle\prod\limits_{l=0}^{\frac{p-1}{2}}}
%EndExpansion
a_{1+2i_{2l+1}}\\
& \\
&  =\chi_{n-1,p+1}\left(  \mathfrak{a}_{0,n-2}\right)  \text{.}%
\end{align*}

Adding all the possible values of $\chi_{n+1,2j+1}\left(  \mathfrak{a}%
_{0,n}\right)  t^{2j+1}$ and putting $t=1$ we get%
\[
\Omega_{n+1}\left(  1,\mathfrak{a}_{0,n}\right)  =a_{n}\Omega_{n}\left(
1,\mathfrak{a}_{0,n-1}\right)  +\Omega_{n-1}\left(  1,\mathfrak{a}%
_{0,n-2}\right)  ,
\]
as desired.
\end{pf}

Obviously, the result on general non-autonomous difference equations can be
further developed for the general non-periodic case. In the present work its
main use is on the study of periodic systems. Anyway, we present here a very
simple example for the general non-autonomous case.

\begin{exmp}
Consider the sequence of Fibonacci matrices $\mathcal{A}=\left\{
a_{i}\right\}  _{i\in\mathbb{%
%TCIMACRO{\U{2115} }%
%BeginExpansion
\mathbb{N}
%EndExpansion
}}$ with $a_{i}=0$ when $i$ is odd. For $n$ odd we have
\[
C_{n}=\left[
\begin{array}
[c]{cc}%
\Omega_{n}\left(  1,\mathfrak{a}_{0,n-1}\right)  & \Omega_{n-1}\left(
1,\mathfrak{a}_{1,n-1}\right) \\
\Omega_{n-1}\left(  1,\mathfrak{a}_{0,n-2}\right)  & \Omega_{n-2}\left(
1,\mathfrak{a}_{1,n-2}\right)
\end{array}
\right]  =\left[
\begin{array}
[c]{cc}%
%TCIMACRO{\dsum \limits_{j=0}^{\frac{n-1}{2}}}%
%BeginExpansion
{\displaystyle\sum\limits_{j=0}^{\frac{n-1}{2}}}
%EndExpansion
a_{2j} & 1\\
1 & 0
\end{array}
\right]  ,
\]
where $\Omega_{n-2}\left(  1,\mathfrak{a}_{1,n-2}\right)  =%
%TCIMACRO{\dsum \limits_{j=0}^{\frac{n-3}{2}}}%
%BeginExpansion
{\displaystyle\sum\limits_{j=0}^{\frac{n-3}{2}}}
%EndExpansion
a_{1+2j}=0$. For $n$ even we have%
\[
C_{n}=\left[
\begin{array}
[c]{cc}%
\Omega_{n}\left(  1,\mathfrak{a}_{0,n-1}\right)  & \Omega_{n-1}\left(
1,\mathfrak{a}_{1,n-1}\right) \\
\Omega_{n-1}\left(  1,\mathfrak{a}_{0,n-2}\right)  & \Omega_{n-2}\left(
1,\mathfrak{a}_{1,n-2}\right)
\end{array}
\right]  =\left[
\begin{array}
[c]{cc}%
1 & 0\\%
%TCIMACRO{\dsum \limits_{j=0}^{\frac{n-2}{2}}}%
%BeginExpansion
{\displaystyle\sum\limits_{j=0}^{\frac{n-2}{2}}}
%EndExpansion
a_{2j} & 1
\end{array}
\right]  .
\]
In this simple example the asymptotic behavior of the solutions depends only
on the convergence of the series $%
%TCIMACRO{\dsum \limits_{j=0}^{\infty}}%
%BeginExpansion
{\displaystyle\sum\limits_{j=0}^{\infty}}
%EndExpansion
a_{2j}$.
\end{exmp}

\section{Periodic Fibonacci difference equation\label{MM}}

In this section we study the periodic generalized Fibonacci difference
equation \ref{EqMForm} using the monodromy matrix and its Floquet multipliers.

\begin{defn}
Consider the periodic non-zero sequence of complex numbers $\mathcal{A}%
=\left\{  a_{i}\right\}  _{i\in\mathbb{%
%TCIMACRO{\U{2115} }%
%BeginExpansion
\mathbb{N}
%EndExpansion
}}$ such that$\mathcal{\ }a_{i+k}=a_{i}$ for any $i\in\mathbb{%
%TCIMACRO{\U{2115} }%
%BeginExpansion
\mathbb{N}
%EndExpansion
}$, some fixed $k\in%
%TCIMACRO{\U{2124} }%
%BeginExpansion
\mathbb{Z}
%EndExpansion
^{+}$ and a periodic sequence of $2\times2$\ Fibonacci matrices%
\[
A_{i}=\left[
\begin{array}
[c]{cc}%
a_{i} & 1\\
1 & 0
\end{array}
\right]  .
\]
We define the monodromy matrix $C_{k}\in%
%TCIMACRO{\U{2102} }%
%BeginExpansion
\mathbb{C}
%EndExpansion
^{2\times2}$ by%
\[
C_{k}\triangleq\prod\limits_{i=0}^{k-1}A_{k-1-i}\text{.}%
\]

\end{defn}

Theorem \ref{Main} naturally applies to the monodromy matrix. The next step is
to determine the Floquet multipliers, which are the eigenvalues of $C_{k}$.
For that purpose we need to compute the determinant and trace of $C_{k}$.

Since $C_{k}=\prod\limits_{i=0}^{k-1}\left[
\begin{array}
[c]{cc}%
a_{k-1-i} & 1\\
1 & 0
\end{array}
\right]  \medskip$, it follows that $\det C_{k}=\left(  -1\right)  ^{k}$. By
theorem \ref{Main} the trace\footnote{The constant $A$ in \cite{EY2011} is
basically $\left(  -1\right)  ^{k}\operatorname{tr}C_{k}$ which is very
dificult to compute without Theorem \ref{Main}.} of the monodromy matrix
$C_{k}$ is%
\[
\operatorname{tr}C_{k}=\Omega_{k}\left(  1,\mathfrak{a}_{0,k-1}\right)
+\Omega_{k-2}\left(  1,\mathfrak{a}_{1,k-2}\right)  .
\]
To understand better this trace we use again the formal indeterminate $t$.

\begin{defn}
The formal polynomial $T_{k}\left(  t\right)  \in%
%TCIMACRO{\U{2102} }%
%BeginExpansion
\mathbb{C}
%EndExpansion
\left[  t\right]  $ is given by%
\[
T_{k}\left(  t\right)  \triangleq\Omega_{k}\left(  t,\mathfrak{a}%
_{0,k-1}\right)  +\Omega_{k-2}\left(  t,\mathfrak{a}_{1,k-2}\right)  .
\]

\end{defn}

$T_{k}\left(  t\right)  $\ is defined such that the trace of $C_{k}$ is
$T_{k}\left(  1\right)  $.

\begin{notation}
We write the polynomial $T_{k}\left(  t\right)  $ as%
\[
T_{k}\left(  t\right)  =%
%TCIMACRO{\dsum \limits_{j=0}^{\frac{k}{2}}}%
%BeginExpansion
{\displaystyle\sum\limits_{j=0}^{\frac{k}{2}}}
%EndExpansion
\Psi_{2j}\left(  \mathfrak{a}_{0,k-1}\right)  t^{2j},
\]
when $k$ is even and%
\[
T_{k}\left(  t\right)  =%
%TCIMACRO{\dsum \limits_{j=0}^{\frac{k-1}{2}}}%
%BeginExpansion
{\displaystyle\sum\limits_{j=0}^{\frac{k-1}{2}}}
%EndExpansion
\Psi_{2j+1}\left(  \mathfrak{a}_{0,k-1}\right)  t^{2j+1},
\]
when $k$ is odd.
\end{notation}

Consider the case of even $k$ (similar reasonings hold for the odd case), we
have%
\begin{align*}
T_{k}\left(  t\right)   &  =%
%TCIMACRO{\dsum \limits_{j=0}^{\frac{k}{2}}}%
%BeginExpansion
{\displaystyle\sum\limits_{j=0}^{\frac{k}{2}}}
%EndExpansion
\chi_{k,2j}\left(  \mathfrak{a}_{0,k-1}\right)  t^{2j}+%
%TCIMACRO{\dsum \limits_{j=0}^{\frac{k}{2}-1}}%
%BeginExpansion
{\displaystyle\sum\limits_{j=0}^{\frac{k}{2}-1}}
%EndExpansion
\chi_{k-2,2j}\left(  \mathfrak{a}_{1,k-2}\right)  t^{2j}\\
& \\
&  =%
%TCIMACRO{\dsum \limits_{j=0}^{\frac{k}{2}}}%
%BeginExpansion
{\displaystyle\sum\limits_{j=0}^{\frac{k}{2}}}
%EndExpansion
\left(  \chi_{k,2j}\left(  \mathfrak{a}_{0,k-1}\right)  +\chi_{k-2,2j}\left(
\mathfrak{a}_{1,k-2}\right)  \right)  t^{2j}\text{,}%
\end{align*}
\bigskip subject to the convention that%
\[
\chi_{k-2,2j}\left(  v\right)  =0\,\text{, if }2j>k-2\text{.}%
\]
Now we write%
\[
T_{k}\left(  t\right)  =%
%TCIMACRO{\dsum \limits_{j=0}^{\frac{k}{2}}}%
%BeginExpansion
{\displaystyle\sum\limits_{j=0}^{\frac{k}{2}}}
%EndExpansion
\Psi_{2j}\left(  \mathfrak{a}_{0,k-1}\right)  t^{2j},
\]
where $\Psi_{2j}\left(  \mathfrak{a}_{0,k-1}\right)  $ is the coefficient of
$t^{2j}$ in $T_{k}\left(  t\right)  $ such that
\[
\Psi_{2j}\left(  \mathfrak{a}_{0,k-1}\right)  =\chi_{k,2j}\left(
\mathfrak{a}_{0,k-1}\right)  +\chi_{k-2,2j}\left(  \mathfrak{a}_{1,k-2}%
\right)  \text{.}%
\]
Please note that the coefficient of the highest degree ($t^{k}$) is%
\[
\Psi_{k}\left(  \mathfrak{a}_{0,k-1}\right)  =\chi_{k,k}\left(  \mathfrak{a}%
_{0,k-1}\right)  +\underset{0}{\underbrace{\chi_{k-2,k}\left(  \mathfrak{a}%
_{1,k-2}\right)  }}=%
%TCIMACRO{\dprod _{i=0}^{k-1}}%
%BeginExpansion
{\displaystyle\prod_{i=0}^{k-1}}
%EndExpansion
a_{i}.
\]

We focus our attention on the monomials in the formal indeterminate $t$.
Remembering that $k$ is even and now $2j<k$. Hence, the coefficient of
$t^{2j}$ is%
\[
\Psi_{2j}\left(  \mathfrak{a}_{0,k-1}\right)  =
\]%
\[
\widehat{\sum\limits_{i_{0},\ldots,i_{2j-1}}^{k,2j}}%
%TCIMACRO{\dprod _{l=0}^{j-1}}%
%BeginExpansion
{\displaystyle\prod_{l=0}^{j-1}}
%EndExpansion
a_{2i_{2l}}%
%TCIMACRO{\dprod _{l=0}^{j-1}}%
%BeginExpansion
{\displaystyle\prod_{l=0}^{j-1}}
%EndExpansion
a_{1+2i_{2l+1}}+\widehat{\sum\limits_{i_{0},\ldots,i_{2j-1}}^{k-2,2j}}%
%TCIMACRO{\dprod _{l=0}^{j-1}}%
%BeginExpansion
{\displaystyle\prod_{l=0}^{j-1}}
%EndExpansion
a_{1+2i_{2l}}%
%TCIMACRO{\dprod _{l=0}^{j-1}}%
%BeginExpansion
{\displaystyle\prod_{l=0}^{j-1}}
%EndExpansion
a_{2+2i_{2l+1}}.
\]
We realize that this rather long expression is nothing but the sum of all
possible products with $j$ factors of coefficients $a_{j}$ such that the first
factor can have index even or odd, the second factor has index odd or even,
respectively, alternating always the parity of the indices along the product
of coefficients. For instance if $k\geq6$ we have in the coefficient of
$t^{4}$ summands which are products of the form $a_{0}a_{1}a_{2}a_{3}$,
$a_{0}a_{1}a_{4}a_{5}$, $a_{0}a_{3}a_{4}a_{5}$, $a_{2}a_{3}a_{4}a_{5}$ or
$a_{1}a_{2}a_{3}a_{4}$, but we do not have the forbidden products $a_{0}%
a_{2}a_{3}a_{4}$ or $a_{1}a_{2}a_{3}a_{5}$, where we would have two even
consecutive indices, in the first case, and odd, in the second case.

We can obtain explicitly the Floquet multipliers in the following result.

\begin{prop}
\label{Floquet}The Floquet multipliers of $C_{k}$ are given by%
\[
\Phi_{k}^{\pm}=\frac{1}{2}\left(  T_{k}(1)\pm\sqrt{\left(  T_{k}(1)\right)
^{2}-4\left(  -1\right)  ^{k}}\right)  .
\]

\end{prop}

\begin{pf}
Since the eigenvalues of any $2\times2$ matrix $A\in%
%TCIMACRO{\U{2102} }%
%BeginExpansion
\mathbb{C}
%EndExpansion
^{2\times2}$ are given by
\[
\lambda_{1,2}=\frac{\operatorname{tr}A\pm\sqrt{\left(  \operatorname{tr}%
A\right)  ^{2}-4\det A}}{2},
\]
it follows that
\[
\Phi_{k}^{\pm}=\frac{T_{k}(1)\pm\sqrt{\left(  T_{k}(1)\right)  ^{2}-4\left(
-1\right)  ^{k}}}{2}.
\]

\end{pf}

Now let $C_{k}=J\Lambda J^{-1}$ be the Jordan canonical form of the monodromy
matrix $C_{k}$. Assume that when $k$ is even we have $|\operatorname{tr}%
C_{k}|\neq2$. Hence, we can write, without loss of generality, the matrices
$J$, $J^{-1}$ and $\Lambda$ as
\[
J=\left[
\begin{array}
[c]{cc}%
J_{11} & J_{12}\\
J_{21} & J_{22}%
\end{array}
\right]  ,\text{ }J^{-1}=\left[
\begin{array}
[c]{cc}%
J_{11}^{-1} & J_{12}^{-1}\\
J_{21}^{-1} & J_{22}^{-1}%
\end{array}
\right]  ,\text{ }\Lambda=\left[
\begin{array}
[c]{cc}%
\Phi_{k}^{-} & 0\\
0 & \Phi_{k}^{+}%
\end{array}
\right]  .
\]

So, we have
\begin{align*}
C_{k}^{n}  &  =J\Lambda^{n} J^{-1}\\
&  = \left[
\begin{array}
[c]{cc}%
J_{11}(\Phi_{k}^{-})^{n}J_{11}^{-1}+ J_{12}(\Phi_{k}^{+})^{n}J_{21}^{-1} &
J_{11}(\Phi_{k}^{-})^{n}J_{12}^{-1}+J_{12}(\Phi_{k}^{+})^{n}J_{22}^{-1}\\
J_{21}(\Phi_{k}^{-})^{n}J_{11}^{-1}+ J_{22}(\Phi_{k}^{+})^{n}J_{21}^{-1} &
J_{21}(\Phi_{k}^{-})^{n}J_{12}^{-1}+J_{22}(\Phi_{k}^{+})^{n}J_{22}^{-1}%
\end{array}
\right]  .
\end{align*}
If $k$ is even and $|\operatorname{tr}C_{k}|=2$, then the eigenvalues of
$C_{k}$ are either $1$ or $-1$ with algebraic multiplicity $2$. In this case
the matrix $\Lambda^{n}$ is either
\[
\left[
\begin{array}
[c]{cc}%
1 & n\\
0 & 1
\end{array}
\right]  \text{ or }\left[
\begin{array}
[c]{cc}%
(-1)^{n} & n(-1)^{n-1}\\
0 & (-1)^{n}%
\end{array}
\right]  .
\]

In the following corollary of Theorem \ref{Main} and Proposition
\ref{Floquet}\ we write explicitly the solution of equation (\ref{naFre}).

\begin{cor}
\label{binet} The solution of the $k$-periodic generalized Fibonacci
difference equation (\ref{naFre}) is given by

\begin{enumerate}
\item Case $k$ odd or ($k$ even and $|\operatorname{tr}C_{k}|\neq2$):%

\begin{align*}
x_{nk}  &  =J_{21}(J_{11}^{-1} x_{1}+J_{21}^{-1}x_{0})(\Phi_{k}^{-}%
)^{n}+J_{22}(J_{21}^{-1}x_{1}+J^{-1}_{22}x_{0})(\Phi_{k}^{+})^{n},\\
x_{nk+1}  &  =J_{11}(J_{11}^{-1} x_{1}+J_{21}^{-1} x_{0})(\Phi_{k}^{-}%
)^{n}+J_{12}(J^{-1}_{21}x_{1}+J^{-1}_{22}x_{0})(\Phi_{k}^{+})^{n},
\end{align*}
and
\[
x_{nk+(i+2)}=a_{i}x_{nk+(i+1)}+x_{nk+i}, i\in\left\{  0,1,2,\ldots
,k-3\right\}  ,
\]
for all $n=0,1,2,\ldots$.

\item Case $k$ even and $\operatorname{tr}C_{k}=2$:%

\begin{align*}
x_{nk}  &  =\left(  J_{21} J^{-1}_{11}+\left(  n J_{21}+J_{22}\right)
J^{-1}_{21}\right)  x_{1}+ \left(  J_{21} J^{-1}_{12}+\left(  n J_{21}%
+J_{22}\right)  J^{-1}_{22}\right)  x_{0},\\
x_{nk+1}  &  =\left(  J_{11} J^{-1}_{11}+\left(  n J_{11}+J_{12}\right)
J^{-1}_{21}\right)  x_{1}+ \left(  J_{11} J^{-1}_{12}+\left(  n J_{11}%
+J_{12}\right)  J^{-1}_{22}\right)  x_{0},
\end{align*}
and
\[
x_{nk+(i+2)}=a_{i}x_{nk+(i+1)}+x_{nk+i}, i\in\left\{  0,1,2,\ldots
,k-3\right\}  ,
\]
for all $n=0,1,2,\ldots$.

\item Case $k$ even and $\operatorname{tr}C_{k}=-2$:%

\begin{align*}
x_{nk}  &  =(-1)^{n}\left(  J_{21} J^{-1}_{11}+\left(  - n J_{21}+
J_{22}\right)  J^{-1}_{21}\right)  x_{1}+\\
&  (-1)^{n}\left(  J_{21} J^{-1}_{12}+\left(  - n J_{21}+ J_{22}\right)
J^{-1}_{22}\right)  x_{0},\\
x_{nk+1}  &  =(-1)^{n}\left(  J_{11} J^{-1}_{11}+\left(  - n J_{11}+
J_{12}\right)  J^{-1}_{21}\right)  x_{1}+\\
&  (-1)^{n}\left(  J_{11} J^{-1}_{12}+\left(  - n J_{11}+ J_{12}\right)
J^{-1}_{22}\right)  x_{0},
\end{align*}
and
\[
x_{nk+(i+2)}=a_{i}x_{nk+(i+1)}+x_{nk+i}, i\in\left\{  0,1,2,\ldots
,k-3\right\}  ,
\]
for all $n=0,1,2,\ldots$.
\end{enumerate}
\end{cor}

\begin{rmk}
Notice that when $x_{0}=0$, $x_{1}=1$ and $k=1$ we have $\Phi_{1}^{+}=\phi$,
$\Phi_{1}^{-}=1-\phi$, $J_{21}=-\phi$, $J_{11}^{-1}=-\frac{1}{\sqrt{5}}%
(1-\phi)$, $J_{22}=1-\phi$ and $J_{21}^{-1}=\frac{1}{\sqrt{5}}\phi$ leading
to
\[
x_{n}=\frac{1}{\sqrt{5}}\left(  \phi^{n}-\left(  1- \phi\right)  ^{n}\right)
, n=0,1,2,\ldots.
\]

\end{rmk}

\section{The structure of the solution\label{Solutions}}

In this section we study the solutions of equation (\ref{naFre}), namely we
determine explicitly the conditions for the periodicity of the solutions in
the case of $k$-periodic Fibonacci equations. The sequence of quotients of
consecutive iterates, such that $q_{0}=\frac{x_{1}}{x_{0}},$ $\ldots,$
$q_{n}=\frac{x_{n+1}}{x_{n}}$ approaches a periodic cycle with period $P$,
where $P$ is a multiple of $k$ (it depends on the Floquet multipliers), which
is related to the convergence to the golden ratio of the quotients of
consecutive iterates in the classic autonomous Fibonacci equation.

\subsection{Odd period}

From Corollary \ref{binet} it follows that
\[
\frac{x_{nk+1}}{x_{nk}}=\frac{\alpha(\Phi_{k}^{-})^{n}+ \beta\left(  \Phi
_{k}^{+}\right)  ^{n}} {\delta(\Phi_{k}^{-})^{n}+\gamma(\Phi_{k}^{+})^{n}},
\]
where $\alpha=J_{11}(J_{11}^{-1} x_{1}+J_{21}^{-1} x_{0})$, $\beta
=J_{12}(J^{-1}_{21}x_{1}+J^{-1}_{22}x_{0})$, $\delta=J_{21}(J_{11}^{-1}
x_{1}+J_{21}^{-1}x_{0})$ and $\gamma=J_{22}(J_{21}^{-1}x_{1}+J^{-1}_{22}%
x_{0})$.

Let the period of equation (\ref{naFre}) be odd, i.e., $k$ is an odd number.
Then
\[
\Phi_{k}^{\pm}=\frac{\operatorname{tr}C_{k}\pm\sqrt{\left(  \operatorname{tr}%
C_{k}\right)  ^{2}+4}}{2}.
\]
Let us assume first that $\operatorname{tr}C_{k}\neq0$. Hence, either
$\left\vert \frac{\Phi_{k}^{-}}{\Phi_{k}^{+}}\right\vert <1$ or $\left\vert
\frac{\Phi_{k}^{-}}{\Phi_{k}^{+}}\right\vert >1$.

If $\left\vert \frac{\Phi_{k}^{-}}{\Phi_{k}^{+}}\right\vert =\left\vert
r\right\vert <1$, we have
\[
\frac{x_{nk+1}}{x_{nk}}=\frac{\alpha r^{n}+\beta}{\delta r^{n}+\gamma
}\longrightarrow\frac{\beta}{\gamma}.
\]
Let $\frac{\beta}{\gamma}=L_{0}$. Similarly, one can show that $\frac
{x_{nk+2}}{x_{nk+1}}\longrightarrow\frac{a_{0}L_{0}+1}{L_{0}}=L_{1}$.
Analogously, $\frac{x_{nk+3}}{x_{nk+2}}\longrightarrow\frac{a_{1}L_{1}%
+1}{L_{1}}=L_{2}$. More generally, $\frac{x_{nk+(i+2)}}{x_{nk+(i+1)}%
}\longrightarrow\frac{a_{i}L_{i}+1}{L_{i}}=L_{i+1}$, $i\in\{0,1,2,\dots
,k-2\}$. Notice that $L_{k+i}=L_{i}$, for all $i=0,1,2,\ldots$. Hence, one can
consider the following $k$-periodic cycle as the limit of the quotients of the
solution of equation (\ref{naFre})
\begin{equation}
\{L_{0},L_{1},L_{2},\ldots,L_{k-1}\}. \label{cyq}%
\end{equation}

\begin{rmk}
Notice that, if there exists an $i\in\{0,1,2,\ldots,k-1\}$ such that $L_{i}%
=0$, then the $k-$periodic cycle (\ref{cyq}) is unbounded since $L_{i+1}%
\rightarrow\infty$. \label{unbc}
\end{rmk}

On the other hand, if $\left\vert \frac{\Phi_{k}^{-}}{\Phi_{k}^{+}}\right\vert
>1$ one can show that $\frac{x_{nk+1}}{x_{nk}}\longrightarrow\frac{\alpha
}{\delta}=\tilde{L}_{0}$, $\frac{x_{nk+2}}{x_{nk+1}}\longrightarrow\frac
{a_{0}\tilde{L}_{0}+1}{\tilde{L}_{0}}=\tilde{L}_{1}$ and more generally
$\frac{x_{nk+(i+2)}}{x_{nk+(i+1)}}\longrightarrow\frac{a_{i}\tilde{L}_{i}%
+1}{\tilde{L}_{i}}=\tilde{L}_{i+1}$, $i\in\{0,1,2,\dots,k-2\}$, yielding the
following $k-$periodic cycle as the limit of the quotients of the solutions
\[
\{\tilde{L}_{0},\tilde{L}_{1},\tilde{L}_{2},\dots,\tilde{L}_{k-1}\}.
\]

Secondly, if $\operatorname{tr}C_{k}=0$, then the eigenvalues of $C_{k}$ are
$1$ and $-1$. Hence, the solutions of equation (\ref{naFre}) are given by
\begin{align*}
x_{kn}  &  =J_{21}(J_{11}^{-1}x_{1}+J_{21}^{-1}x_{0})(-1)^{n}+J_{22}%
(J_{21}^{-1}x_{1}+J_{22}^{-1}x_{0}),\\
x_{kn+1}  &  =J_{11}(J_{11}^{-1}x_{1}+J_{21}^{-1}x_{0})(-1)^{n}+J_{12}%
(J_{21}^{-1}x_{1}+J_{22}^{-1}x_{0})
\end{align*}
and
\[
x_{nk+(i+2)}=a_{i}x_{nk+(i+1)}+x_{nk+i},i\in\left\{  0,1,2,\ldots,k-3\right\}
,
\]
for all $n=0,1,2,\ldots$.

A simple computation shows that $\{x_{0},x_{1},x_{2},\ldots,x_{2k-1}\}$ is a
$2k$-periodic solution of equation (\ref{naFre}), where
\begin{align*}
x_{0}  &  =J_{21}(J_{11}^{-1}x_{1}+J_{21}^{-1}x_{0})+J_{22}(J_{21}^{-1}%
x_{1}+J_{22}^{-1}x_{0}),\\
x_{1}  &  =-J_{11}(J_{11}^{-1}x_{1}+J_{21}^{-1}x_{0})+J_{12}(J_{21}^{-1}%
x_{1}+J_{22}^{-1}x_{0}),\\
x_{i+2}  &  =a_{i}x_{i+1}+x_{i},i\in\left\{  0,1,2,\ldots,k-3\right\}  ,\\
x_{k}  &  =-J_{21}(J_{11}^{-1}x_{1}+J_{21}^{-1}x_{0})+J_{22}(J_{21}^{-1}%
x_{1}+J_{22}^{-1}x_{0}),\\
x_{k+1}  &  =J_{11}(J_{11}^{-1}x_{1}+J_{21}^{-1}x_{0})+J_{12}(J_{21}^{-1}%
x_{1}+J_{22}^{-1}x_{0}),\\
x_{k+(i+2)}  &  =a_{i}x_{k+(i+1)}+x_{k+i},i\in\left\{  0,1,2,\ldots
,k-3\right\}  .
\end{align*}

\subsection{Even period\label{speven}}

Let the period of equation (\ref{naFre}) be even. First let us assume that
$|\operatorname{tr}C_{k}|<2$. Then
\[
\Phi_{k}^{\pm}=\frac{\operatorname{tr}C_{k}\pm\Im\sqrt{4-\left(
\operatorname{tr}C_{k}\right)  ^{2}}}{2},
\]
with $|\Phi_{k}^{+}|=|\Phi_{k}^{-}|=1$ and $\Im$ being the imaginary unit.
Hence, the eigenvalues of $C_{k}$ lie on the unit circle. Moreover, the matrix
$\Lambda$ is periodic being
\[
\frac{\operatorname{tr}C_{k}\pm\Im\sqrt{4-\left(  \operatorname{tr}%
C_{k}\right)  ^{2}}}{2}=e^{\Im\theta},
\]
with $\theta=\arctan\frac{\sqrt{4-\left(  \operatorname{tr}C_{k}\right)  ^{2}%
}}{\operatorname{tr}C_{k}}$, $\theta\in\lbrack-\frac{\pi}{2},\frac{\pi}{2}]$.
Let $P=\frac{2\pi}{\theta}$ be the period of the monodromy matrix, i.e.,
$C_{k}^{P}=C_{k}$. Hence, the solutions of equation (\ref{naFre}) are
periodic. Let us now determine this period.

From Corollary \ref{binet} it follows that
\begin{align*}
x_{0}  &  =J_{21}(J_{11}^{-1}x_{1}+J_{21}^{-1}x_{0})+J_{22}(J_{21}^{-1}%
x_{1}+J_{22}^{-1}x_{0}),\\
x_{1}  &  =J_{11}(J_{11}^{-1}x_{1}+J_{21}^{-1}x_{0})+J_{12}(J_{21}^{-1}%
x_{1}+J_{22}^{-1}x_{0}),\\
x_{i+2}  &  =a_{i}x_{i+1}+x_{i},i\in\left\{  0,1,2,\ldots,k-3\right\}  ,\\
x_{k}  &  =J_{21}(J_{11}^{-1}x_{1}+J_{21}^{-1}x_{0})\Phi_{k}^{-}+J_{22}%
(J_{21}^{-1}x_{1}+J_{22}^{-1}x_{0})\Phi_{k}^{+},\\
x_{k+1}  &  =J_{11}(J_{11}^{-1}x_{1}+J_{21}^{-1}x_{0})\Phi_{k}^{-}%
+J_{12}(J_{21}^{-1}x_{1}+J_{22}^{-1}x_{0})\Phi_{k}^{+},\\
x_{k+(i+2)}  &  =a_{i}x_{k+(i+1)}+x_{k+i},i\in\left\{  0,1,2,\ldots
,k-3\right\}  ,\\
&  \ldots\\
x_{(P-1)k}  &  =J_{21}(J_{11}^{-1}x_{1}+J_{21}^{-1}x_{0})(\Phi_{k}^{-}%
)^{P-1}+J_{22}(J_{21}^{-1}x_{1}+J_{22}^{-1}x_{0})(\Phi_{k}^{+})^{P-1},\\
x_{(P-1)k+1}  &  =J_{11}(J_{11}^{-1}x_{1}+J_{21}^{-1}x_{0})(\Phi_{k}%
^{-})^{P-1}+J_{12}(J_{21}^{-1}x_{1}+J_{22}^{-1}x_{0})(\Phi_{k}^{+})^{P-1},\\
x_{(P-1)k+(i+2)}  &  =a_{i}x_{(P-1)k+(i+1)}+x_{(P-1)k+i},i\in\left\{
0,1,2,\ldots,k-3\right\}  .
\end{align*}
Since the monodromy matrix is $P$-periodic, it follows that $(\Phi_{k}^{\pm
})^{P+i}=(\Phi_{k}^{\pm})^{i}$, for all $i=0,1,2,\ldots$. Hence,
$x_{kP+i}=x_{i}$, for all $i=0,1,2,\ldots$. This implies that the minimal
period of the cycle is $kP$.

In conclusion, if $|\operatorname{tr}C_{k}|<2$, the following cycle is a
$kP$-periodic solution of equation (\ref{naFre})
\[
\{x_{0},x_{1},x_{2},\ldots,x_{kP-1}\}.
\]

Second, if $|\operatorname{tr}C_{k}|=2$, it follows that $\left\vert
\frac{\Phi_{k}^{-}}{\Phi_{k}^{+}}\right\vert =1$. In this case we will have a
$k-$periodic cycle as the limit of the quotients of the solution.

Third, if $|\operatorname{tr}C_{k}|>2$, it follows that we will have either
$\left\vert \frac{\Phi_{k}^{-}}{\Phi_{k}^{+}}\right\vert <1$ or $\left\vert
\frac{\Phi_{k}^{-}}{\Phi_{k}^{+}}\right\vert >1$. Following the same ideas as
in the odd case we will have a $k$-periodic cycle as the limit of the
quotients of the solution of equation (\ref{naFre}).

\section{Applications and examples \label{exampls}}

In this section we study some examples of the periodic generalized Fibonacci
difference equation (\ref{naFre}). We start by a $3$-periodic equation.

\begin{exmp}
Let the period of equation (\ref{naFre}) be 3, i.e., $k=3$. The monodromy
matrix is given by
\[
C_{3}=\left[
\begin{array}
[c]{cc}%
a_{0}+a_{2}+a_{0}a_{1}a_{2} & 1+a_{1}a_{2}\\
1+a_{0}a_{1} & a_{1}%
\end{array}
\right]  ,
\]
yielding the following representation%
\[
C_{3}^{n}=J\Lambda^{n}J^{-1}%
\]
with
\begin{align*}
J  &  =\left[
\begin{array}
[c]{cc}%
\frac{\Delta-2a_{1}-\sqrt{4+\Delta^{2}}}{2(1+a_{0}a_{1})} & \frac
{\Delta-2a_{1}+\sqrt{4+\Delta^{2}}}{2(1+a_{0}a_{1})}\\
1 & 1
\end{array}
\right]  ,\\
J^{-1}  &  =\left[
\begin{array}
[c]{cc}%
-\frac{1+a_{0}a_{1}}{\sqrt{4+\Delta^{2}}} & \frac{\Delta-2a_{1}+\sqrt
{4+\Delta^{2}}}{2\sqrt{4+\Delta^{2}}}\\
\frac{1+a_{0}a_{1}}{\sqrt{4+\Delta^{2}}} & -\frac{\Delta-2a_{1}-\sqrt
{4+\Delta^{2}}}{2\sqrt{4+\Delta^{2}}}%
\end{array}
\right]  ,\\
\Lambda^{n}  &  =\left[
\begin{array}
[c]{cc}%
\left(  \frac{\Delta-\sqrt{4+\Delta^{2}}}{2}\right)  ^{n} & 0\\
0 & \left(  \frac{\Delta+\sqrt{4+\Delta^{2}}}{2}\right)  ^{n}%
\end{array}
\right]  ,
\end{align*}
where $\Delta=a_{0}+a_{1}+a_{2}+a_{0}a_{1}a_{2}$.

\begin{enumerate}
\item Case $\Delta=\operatorname{tr}C_{3}=0$. This occurs when $a_{0}%
=-\frac{a_{1}+a_{2}}{1+a_{1}a_{2}}$. Hence, the general solution is given by
\[
x_{3n}=\frac{1+(-1)^{n}+(1-(-1)^{n})a_{1}}{2}x_{0}+\frac{(1-(-1)^{n}%
)(1-a_{1}^{2})}{2(1+a_{1}a_{2})}x_{1},
\]%
\[
x_{3n+1}=\frac{1-(-1)^{n}}{2}(1+a_{1}a_{2})x_{0}+\frac{1+(-1)^{n}%
-(1-(-1)^{n})a_{1}}{2}x_{1}%
\]
and
\[
x_{3n+2}=a_{0}x_{3n+1}+x_{3n},
\]
for all $n=0,1,2,\ldots$. It is straightforward to see that this general
solution is, in fact, a $6$-periodic solution of the form%

%TCIMACRO{\FRAME{ftFU}{2.9776in}{1.8887in}{0pt}{\Qcb{A $6-$periodic solution
%($\{5,1,4,5,-6,11\}$) of a $3-$periodic equation. The values of the parameters
%are $a_{0}=-1$, $a_{1}=1$ and $a_{2}=-2$ with initial conditions $x_{0}=5$ and
%$x_{1}=1$. In this case the eigenvalues of $C_{3}$ are $1$ and $-1$ since
%$\operatorname{tr}C_{k}=0$.}}{\Qlb{fig:orbiteven_p}}{orbit_even_3p1.eps}%
%maintain-aspect-ratio TRUE;  display "FULL";  valid_file "F";
%width 2.9776in;  height 1.8887in;  depth 0pt;  original-width 3.5552in;
%original-height 2.2453in;  cropleft "0";  croptop "1";  cropright "1";
%cropbottom "0";  filename 'orbit_even_3p1.eps';file-properties "XNPEU";}}}%
%BeginExpansion
\begin{figure}
[t]
\begin{center}
\includegraphics[
height=1.8887in,
width=2.9776in
]%
{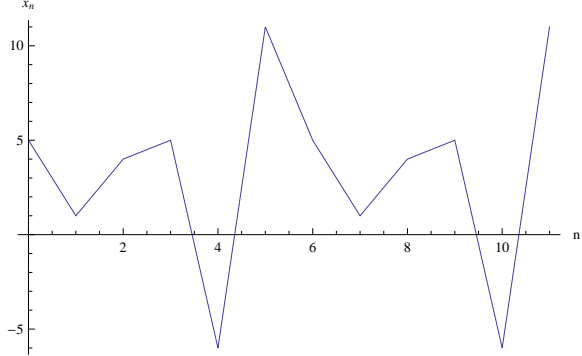}%
\caption{A $6-$periodic solution ($\{5,1,4,5,-6,11\}$) of a $3-$periodic
equation. The values of the parameters are $a_{0}=-1$, $a_{1}=1$ and
$a_{2}=-2$ with initial conditions $x_{0}=5$ and $x_{1}=1$. In this case the
eigenvalues of $C_{3}$ are $1$ and $-1$ since $\operatorname{tr}C_{k}=0$.}%
\label{fig:orbiteven_p}%
\end{center}
\end{figure}
%EndExpansion%
\begin{align*}
&  \{x_{0},x_{1},a_{0}x_{1}+x_{0},a_{1}x_{0}+\frac{1-a_{1}^{2}}{2(1+a_{1}%
a_{2})}x_{1},(1+a_{1}a_{2})x_{0}-a_{1}x_{1},\\
&  a_{0}((1+a_{1}a_{2})x_{0}-a_{1}x_{1})+a_{1}x_{0}+\frac{1-a_{1}^{2}%
}{2(1+a_{1}a_{2})}x_{1}\}.
\end{align*}
Notice that this cycle can be not bounded if $a_{1}a_{2}=-1$. In Figure
\ref{fig:orbiteven_p} we can see an example of this situation. Using a
$3$-periodic sequence of parameters we present in the plane $(n,x_{n})$ a
$6$-periodic solution of equation (\ref{naFre}).

\item Case $\Delta=\operatorname{tr}C_{3}\neq0$. The general solution of the
$3$-periodic equation is given by%
\[
^{x_{3n}=\frac{1+a_{0}a_{1}}{\sqrt{4+\Delta^{2}}}\left(  \left(  \Phi_{3}%
^{+}\right)  ^{n}-\left(  \Phi_{3}^{-}\right)  ^{n}\right)  x_{1}+\left(
\frac{\Delta-2a_{1}+\sqrt{4+\Delta^{2}}}{2\sqrt{4+\Delta^{2}}}\left(  \Phi
_{3}^{-}\right)  ^{n}+\frac{-\Delta+2a_{1}+\sqrt{4+\Delta^{2}}}{2\sqrt
{4+\Delta^{2}}}\left(  \Phi_{3}^{+}\right)  ^{n}\right)  x_{0},}%
\]%
\[
^{x_{3n+1}=\frac{1+a_{1}a_{2}}{\sqrt{4+\Delta^{2}}}\left(  \left(  \Phi
_{3}^{+}\right)  ^{n}-\left(  \Phi_{3}^{-}\right)  ^{n}\right)  x_{0}+\left(
\frac{-\Delta+2a_{1}+\sqrt{4+\Delta^{2}}}{2\sqrt{4+\Delta^{2}}}\left(
\Phi_{3}^{-}\right)  ^{n}+\frac{\Delta-2a_{1}+\sqrt{4+\Delta^{2}}}%
{2\sqrt{4+\Delta^{2}}}\left(  \Phi_{3}^{+}\right)  ^{n}\right)  x_{1},}%
\]
and
\[
x_{3n+2}=a_{0}x_{3n+1}+x_{3n},
\]
for all $n=0,1,2,\ldots$, where $\Phi_{3}^{\pm}=\frac{\Delta\pm\sqrt
{4+\Delta^{2}}}{2}$.%

%TCIMACRO{\FRAME{ftFU}{3.0934in}{1.9363in}{0pt}{\Qcb{A $3-$periodic cycle
%$\{-3.61803,0.723607,-0.618034\}$ for the ratios of the solution of a
%$3-$periodic equation. The values of the parameters are $a_{0}=1$, $a_{1}=-2$
%and $a_{2}=-2$ with initial conditions $x_{0}=5$ and $x_{1}=1$.}%
%}{\Qlb{fig:orbiteven}}{orbit_even_3p2.eps}%
%width 3.0934in;  height 1.9363in;  depth 0pt;  original-width 3.5552in;
%original-height 2.2453in;  cropleft "-0.001373";  croptop "1";
%cropright "0.998627";  cropbottom "0";
%filename 'orbit_even_3p2.eps';file-properties "XNPEU";}}}%
%BeginExpansion
\begin{figure}
[t]
\begin{center}
\includegraphics[
trim=-0.004881in 0.000000in 0.004881in 0.000000in,
height=1.9363in,
width=3.0934in
]%
{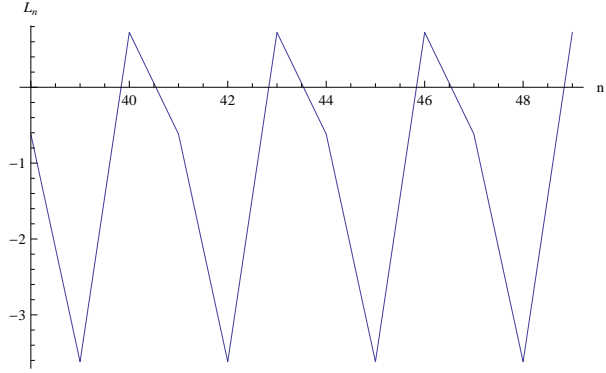}%
\caption{A $3-$periodic cycle $\{-3.61803,0.723607,-0.618034\}$ for the ratios
of the solution of a $3-$periodic equation. The values of the parameters are
$a_{0}=1$, $a_{1}=-2$ and $a_{2}=-2$ with initial conditions $x_{0}=5$ and
$x_{1}=1$.}%
\label{fig:orbiteven}%
\end{center}
\end{figure}
%EndExpansion

\begin{itemize}
\item Case $\Delta>0$. If the trace of the monodromy matrix is positive, then
the limit of the quotients of the solution converge to the following
$3$-periodic cycle
\begin{equation}
\{L_{0},\frac{a_{0}L_{0}+1}{L_{0}},\frac{a_{1}(a_{0}L_{0}+1)+L_{0}}{a_{0}%
L_{0}+1}\}, \label{eq:c1}%
\end{equation}
where
\[
L_{0}=\frac{2(1+a_{1}a_{2})x_{0}+(\Delta-2a_{1}+\sqrt{4+\Delta^{2}})x_{1}%
}{2(1+a_{0}a_{1})x_{1}+(-\Delta+2a_{1}+\sqrt{4+\Delta^{2}})x_{0}}.
\]

\item Finally, $\Delta<0$. In this case we obtain a cycle as in (\ref{eq:c1})
with
\[
L_{0}=\frac{-2(1+a_{1}a_{2})x_{0}+(-\Delta+2a_{1}+\sqrt{4+\Delta^{2}})x_{1}%
}{-2(1+a_{0}a_{1})x_{1}+(\Delta-2a_{1}+\sqrt{4+\Delta^{2}})x_{0}}.
\]

\end{itemize}

In Figure \ref{fig:orbiteven} we have an example of this case. Using a
$3$-periodic sequence of parameters we plot in the plane $(n,\frac{x_{n+1}%
}{x_{n}})$ a $3$-periodic cycle for the quotients of the solutions.
\end{enumerate}
\end{exmp}

\begin{exmp}
Let us now consider a $4$-periodic equation. The monodromy matrix is given by
\[
C_{4}=\left[
\begin{array}
[c]{cc}%
1+a_{0}a_{1}+a_{0}a_{3}+a_{2}a_{3}+a_{0}a_{1}a_{2}a_{3} & a_{1}+a_{3}%
+a_{1}a_{2}a_{3}\\
a_{0}+a_{2}+a_{0}a_{1}a_{2} & 1+a_{1}a_{2}%
\end{array}
\right]  .
\]

\begin{enumerate}
\item Consider the case $\operatorname{tr}C_{4}=2$. This occurs when
$a_{0}=-\frac{a_{2}(a_{1}+a_{3})}{a_{1}+a_{3}+a_{1}a_{2}a_{3}}$. Hence, the
monodromy matrix can be simplified as
\[
C_{4}=\left[
\begin{array}
[c]{cc}%
1-a_{1}a_{2} & a_{1}+a_{3}+a_{1}a_{2}a_{3}\\
-\frac{a_{1}^{2}a_{2}^{2}}{a_{3}+a_{1}\left(  1+a_{2}a_{3}\right)  } &
1+a_{1}a_{2}%
\end{array}
\right]  ,
\]
yielding the following representation%
\[
C_{4}^{n}=\left[
\begin{array}
[c]{cc}%
\frac{a_{1}+a_{3}+a_{1}a_{2}a_{3}}{a_{1}a_{2}} & \frac{-a_{1}-a_{3}-a_{1}%
a_{2}a_{3}}{a_{1}^{2}a_{2}^{2}}\\
1 & 0
\end{array}
\right]  \left[
\begin{array}
[c]{cc}%
1 & n\\
0 & 1
\end{array}
\right]  \left[
\begin{array}
[c]{cc}%
0 & 1\\
-\frac{a_{1}^{2}a_{2}^{2}}{a_{3}+a_{1}\left(  1+a_{2}a_{3}\right)  } &
a_{1}a_{2}%
\end{array}
\right]  .
\]
It follows that
\begin{align*}
x_{4n}  &  =(1+na_{1}a_{2})x_{0}-\frac{a_{1}^{2}a_{2}^{2}n}{a_{3}+a_{1}\left(
1+a_{2}a_{3}\right)  }x_{1},\\
x_{4n+1}  &  =(a_{1}+a_{3}+a_{1}a_{2}a_{3})nx_{0}+(1-a_{1}a_{2}n)x_{1},\\
x_{4n+2}  &  =-\frac{a_{2}(a_{1}+a_{3})}{a_{1}+a_{3}+a_{1}a_{2}a_{3}}%
x_{4n+1}+x_{4n},\\
x_{4n+3}  &  =a_{1}x_{4n+2}+x_{4n+1}.
\end{align*}
With some computations we see that $\frac{x_{4n+1}}{x_{4n}}\rightarrow L_{0}$,
$\frac{x_{4n+2}}{x_{4n+1}}\rightarrow L_{1}=\frac{a_{0}L_{0}+1}{L_{0}}$,
$\frac{x_{4n+3}}{x_{4n+2}}\rightarrow L_{2}=\frac{a_{1}L_{1}+1}{L_{1}}$ and
$\frac{x_{4n+4}}{x_{4n+3}}\rightarrow L_{3}=\frac{a_{2}L_{2}+1}{L_{2}}$,
where
\[
L_{0}=\frac{(a_{1}+a_{3}+a_{1}a_{3})^{2}x_{0}-a_{1}a_{2}(a_{1}+a_{3}%
+a_{1}a_{3})x_{1}}{a_{1}a_{2}(a_{1}+a_{3}+a_{1}a_{3})x_{0}-(a_{1}a_{2}%
)^{2}x_{1}}.
\]
Hence, $\{L_{0},L_{1},L_{2},L_{3}\}$ is a $4$-periodic cycle of the limiting
process of the quotients of the solution of the equation. See Figure
\ref{fig:oddtr=2} for a concrete example.%

%TCIMACRO{\FRAME{ftFU}{3.0511in}{1.8853in}{0pt}{\Qcb{A $4-$periodic cycle of
%the limiting of the ratios of the solution of a $4-$periodic generalized
%Fibonacci equation. In this example we consider $a_{0}=-1$, $a_{1}=1$,
%$a_{2}=-1$ and $a_{3}=-2$ with initial conditions $x_{0}=5$ and $x_{1}=1$
%originating the cycle $\{-1,-2,0.5,1\}$.}}{\Qlb{fig:oddtr=2}}%
%{orbit_odd_tr=2.eps}{\special{ language "Scientific Word";  type "GRAPHIC";
%maintain-aspect-ratio TRUE;  display "FULL";  valid_file "F";
%width 3.0511in;  height 1.8853in;  depth 0pt;  original-width 3.5552in;
%original-height 2.1872in;  cropleft "0";  croptop "1";  cropright "1";
%cropbottom "0";  filename 'orbit_odd_tr=2.eps';file-properties "XNPEU";}}}%
%BeginExpansion
\begin{figure}
[t]
\begin{center}
\includegraphics[
height=1.8853in,
width=3.0511in
]%
{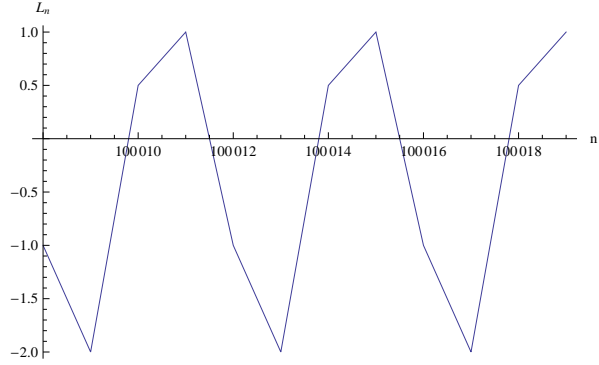}%
\caption{A $4-$periodic cycle of the limiting of the ratios of the solution of
a $4-$periodic generalized Fibonacci equation. In this example we consider
$a_{0}=-1$, $a_{1}=1$, $a_{2}=-1$ and $a_{3}=-2$ with initial conditions
$x_{0}=5$ and $x_{1}=1$ originating the cycle $\{-1,-2,0.5,1\}$.}%
\label{fig:oddtr=2}%
\end{center}
\end{figure}
%EndExpansion

\item Consider $\operatorname{tr}C_{4}=-2$. This occurs when $a_{0}%
=-\frac{4+a_{2}(a_{1}+a_{3})}{a_{1}+a_{3}+a_{1}a_{2}a_{3}}$. The general
solution is given by
\[
x_{4n}=(-1)^{n}\left(  \left(  1-(2+a_{1}a_{2})n\right)  x_{0}+\frac{n\left(
2+a_{1}a_{2}\right)  ^{2}}{a_{3}+a_{1}\left(  1+a_{2}a_{3}\right)  }%
x_{1}\right)  ,
\]%
\[
x_{4n+1}=(-1)^{n}\left(  \left(  1+(2+a_{1}a_{2})n\right)  x_{1}-\left(
a_{1}+a_{3}+a_{1}a_{2}a_{3}\right)  nx_{0}\right)  ,
\]%
\[
x_{4n+2}=a_{0}x_{4n+1}+x_{4n}%
\]
and
\[
x_{4n+3}=a_{1}x_{4n+2}+x_{4n+1},
\]
for all $n=0,1,2,\ldots$. Hence, $\frac{x_{4n+1}}{x_{4n}}\rightarrow L_{0}$,
$\frac{x_{4n+2}}{x_{4n+1}}\rightarrow L_{1}=\frac{a_{0}L_{0}+1}{L_{0}}$,
$\frac{x_{4n+3}}{x_{4n+2}}\rightarrow L_{2}=\frac{a_{1}L_{1}+1}{L_{1}}$ and
$\frac{x_{4n+4}}{x_{4n+3}}\rightarrow L_{3}=\frac{a_{2}L_{2}+1}{L_{2}}$,
where
\[
L_{0}=\frac{(a_{1}+a_{3}+a_{1}a_{2}a_{3})((2+a_{1}a_{2})x_{1}-(a_{1}%
+a_{3}+a_{1}a_{2}a_{3})x_{0})}{(2+a_{1}a_{2})^{2}x_{1}-(2+a_{1}a_{2}%
)(a_{1}+a_{3}+a_{1}a_{2}a_{3})x_{0}}.
\]

With initial conditions $x_{0}=5$ and $x_{1}=1$, for instance, we can find a
$4$-periodic solution $\{10.625,-0.211765,0.277775,1.60003\}$ when
$a_{0}=-\frac{6}{17}$, $a_{1}=5$, $a_{2}=-2$ and $a_{3}=10$ and an unbounded
cycle $\{1,0,\infty,-1\}$ when $a_{0}=-1$, $a_{1}=1$, $a_{2}=-1$ and $a_{3}=2$.

\item Consider $|\operatorname{tr}C_{4}|\neq2$. The general solution of the
$4$-periodic equation is given by
\[
^{%
\begin{array}
[c]{ll}%
x_{4n}= & \frac{2(a_{0}+a_{2}+a_{0}a_{1}a_{2})\left(  \left(  \Phi_{4}%
^{+}\right)  ^{n}-\left(  \Phi_{4}^{-}\right)  ^{n}\right)  }{2\sqrt
{\Delta^{2}-4}}x_{1}+\\
& \frac{(\Delta-2-2a_{1}a_{2}+\sqrt{\Delta^{2}-4})\left(  \Phi_{4}^{-}\right)
^{n}-(\Delta-2-2a_{1}a_{2}-\sqrt{\Delta^{2}-4})\left(  \Phi_{4}^{+}\right)
^{n}}{2\sqrt{\Delta^{2}-4}}x_{0},
\end{array}
}%
\]%
\[
^{%
\begin{array}
[c]{ll}%
x_{4n+1}= & \frac{2(a_{1}+a_{3}+a_{1}a_{2}a_{3})\left(  \left(  \Phi_{4}%
^{+}\right)  ^{n}-\left(  \Phi_{4}^{-}\right)  ^{n}\right)  }{2\sqrt
{\Delta^{2}-4}}x_{0}+\\
& +\frac{(\Delta-2-2a_{1}a_{2}+\sqrt{\Delta^{2}-4})\left(  \Phi_{4}%
^{+}\right)  ^{n}-(\Delta-2-2a_{1}a_{2}-\sqrt{\Delta^{2}-4})\left(  \Phi
_{4}^{-}\right)  ^{n}}{2\sqrt{\Delta^{2}-4}}x_{1},
\end{array}
}%
\]%
\[
x_{4n+2}=a_{0}x_{4n+1}+x_{4n}%
\]
and
\[
x_{4n+3}=a_{1}x_{4n+2}+x_{4n+1},
\]
for all $n=0,1,2,\ldots$, where $\Delta$ is the trace of $C_{4}$, i.e.,
$\Delta=2+a_{0}a_{1}+a_{1}a_{2}+a_{0}a_{3}+a_{2}a_{3}+a_{0}a_{1}a_{2}a_{3}$
and $\Phi_{4}^{\pm}=\frac{\Delta\pm\sqrt{\Delta^{2}-4}}{2}$.

\begin{itemize}
\item Case $|\Delta|>2$. The quotients of the solution converge to a
$4$-periodic cycle of the form
\[
\{L_{0},\frac{a_{0}L_{0}+1}{L_{0}},a_{1}+\frac{L_{0}}{a_{0}L_{0}+1}%
,a_{2}+\frac{a_{0}L_{0}+1}{a_{1}(a_{0}L_{0}+1)+L_{0}}\},
\]
where $L_{0}=\lim\frac{x_{4n+1}}{x_{4n}}$. Notice that this cycle may be unbounded.

\item Case $|\Delta|<2$. Since we have a pair of complex conjugated Floquet
multipliers, the solution of the $4$-periodic equation is periodic with period
$4P$, $P=\frac{2\pi}{\theta}$ where $\theta=\arctan\frac{\sqrt{4-\Delta^{2}}%
}{\Delta}$, $\theta\in\left[  -\frac{\pi}{2},\frac{\pi}{2}\right]  $. See
Subsection \ref{speven} for the computations of this cycle.

Let us present a concrete example to illustrate this case. Let us consider a
$4$-periodic sequence of parameters $a_{n}=a_{n+4}$, for all $n=0,1,2,\ldots$
such that $a_{0}=-1$, $a_{1}=1$, $a_{2}=-1$ and $a_{3}=-1$. With some
computations we see that we have $\theta=\arctan\sqrt{3}=\frac{\pi}{3}$ and
$P=6$. Hence the monodromy matrix $C_{4}$ is $6$-periodic. Consequently, we
have a $24$-periodic cycle as the solution of this equation. A straightforward
computation shows that this cycle is given by%
\begin{align*}
&  \{x_{0},x_{1},x_{0}-x_{1},x_{0},-x_{1},x_{0}+x_{1},-x_{0}-2x_{1},-x_{1},\\
&  -x_{0}-x_{1},x_{0},-2x_{0}-x_{1},-x_{0}-x_{1},-x_{0},-x_{1},-x_{0}%
+x_{1},-x_{0},x_{1},\\
&  -x_{0}-x_{1},x_{0}+2x_{1},x_{1},x_{0}+x_{1},-x_{0},2x_{0}+x_{1},x_{0}%
+x_{1}\}.
\end{align*}
Notice that any cyclic permutation of this last sequence of parameters leads
necessarily to a $24$-periodic solution.
\end{itemize}
\end{enumerate}
\end{exmp}

\begin{akm}
This research was partially supported by the Funda\c{c}\~{a}o para a
Ci\^{e}ncia e Tecnologia FCT/Portugal.

We want to thank the contribution of the two anonymous referees with valuable
remarks and corrections which improved this work in its revised version.
Moreover, the second referee pointed the valuable references \cite{Jesus,
Mallik1997, Mallik1998, Mallik2000} to us.
\end{akm}

%

%TCIMACRO{\TeXButton{Appendix}{\appendix}}%
%BeginExpansion
\appendix
%EndExpansion



%% References


%\bibliographystyle{elsarticle-num}



%\begin{thebibliography}{00}


%\bibitem{}

%\end{thebibliography}

\end{document}